\documentclass[final]{siamltex}

\usepackage{amssymb}
\usepackage{amsmath}
\usepackage{theorem}
\usepackage{algorithm}
\usepackage{graphicx}
\usepackage[usenames]{color}
\usepackage{multirow}
\usepackage[normalem]{ulem}
\usepackage{enumitem}
\usepackage{dsfont}

\newcommand{\R}{\mathbb{R}}
\newcommand{\N}{\mathbb{N}}

\newcommand{\ve}[1]{ #1}
\newtheorem{Thm}{Theorem}[section]
\newtheorem{Def}{Definition}[section]
\newtheorem{Cor}{Corollary}[section]
\newtheorem{Lemma}{Lemma}[section]
\newtheorem{Pro}{Proposition}[section]

  {\end{list}}

\theoremheaderfont{\scshape}
\makeatletter

  \gdef\listctr{list\romannumeral\the\@listdepth}\expandafter
  
\makeatother

  {\end{list}}

\def\proof{{\it Proof. }}
\def\x{ {\ve x}} \def\s{\ve\sigma} \def\sk{\ve \sigma^{(k)}}
\def\vesi{{\s}}
\def\vesik{\sk}

\def\y{ {\ve y}} 
\def\u{ {\ve u}}\def\lamk{{\lambda^{(k)}}}
\def\w{{\ve w}}
\def\z{ {\ve z}}

\def\d{{\ve d}}\def\p{{\ve p}}
\def\xk{ \x^{(k)}}
\def\yk{ \y^{(k)}}

\def\zk{ \z^{(k)}}\def\barx{\bar{\x}}
\def\yk{ \y^{(k)}}
\def\dk{ \d^{(k)}}
\def\xkk{ \x^{(k+1)}}

\def\ak{{\alpha_k}}
\def\ek{{\epsilon_k}}

\def\kinN{{k\in\mathbb N}}

\def\dom{\mathrm{dom}}

\def\R{\mathbb R}

\def\endproof{\hfill$\square$\vspace{0.3cm}\\}
\def\vesi{{\s}}
\def\vesik{\sk}
\def\hs{h_{\vesi}} \def\dhs{\hs'}
\def\ths{\tilde h_{\vesi,\gamma}}  \newcommand{\thsxy}[2]{\ths(#1,#2)}
\def\thsk{\tilde h_{\vesik,\gamma}}\newcommand{\thskxy}[2]{\thsk(#1,#2)} 
\def\hsk{h_{\vesik}} \def\hsbar{h_{\bar\vesi}}
\newcommand{\hsxy}[2]{\hs(#1,#2)} \newcommand{\dhsxy}[2]{\dhs(#1,#1;#2)}
\newcommand{\hskxy}[2]{\hsk(#1,#2)}
\def\vesij{\s^{(j)}}
\def\thsj{\tilde h_{\vesij,\gamma}}\newcommand{\thsjxy}[2]{\thsj(#1,#2)}

\def\ds{d_{\vesi}}
\def\dsk{d_{\vesik}}

\def\kinN{{k\in\N}}
\def\prox{{\mathrm{prox}}}
\def\k{^{(k)}}

\def\tyk{\tilde \y^{(k)}}
\def\tykl{\tilde\y^{(k,l)}}\def\bykl{\bar\y^{(k,l)}}
\def\Fsk{F_{\sk}}\def\Psisk{\Psi_{\sk}}
\newcommand{\Psiskxy}[1]{\Psisk(#1,\xk)}
\def\v{\ve v} \def\vl{\v^{(l)}} 
\def\u{\ve u} 
\def\G{{\mathcal G}}
\newcommand{\Gskxy}[2]{\G_{\sk}(#1,#2,\xk)}
\def\w{\ve w}\def\hatx{\hat\x}\def\Dkk{{D_{k+1}}}\def\Dk{{D_k}}
\def\xjj{\x^{(j+1)}}
\def\ty{\tilde \y}
  
\def\bthsk{\tilde h_{\vesik,\bar\gamma}}  \newcommand{\bthskxy}[2]{\bthsk(#1,#2)}
\def\bR{\bar {\R}}\def\puntino{\ \cdot\ } 

\def\lamj{\lambda^{(j)}} \def\xj{\x^{(j)}} \def\tyj{\ty^{(j)}} 
\makeatletter
\newcommand{\projalg}[1]{\ifcase #1\or LS\or VMILA\or PDHG\or PDHG-D\or CBGP\or SGP\or SS\or WSS2%
\else \@ctrerr \fi}
\makeatother

\newcounter{ilchangectr}

\definecolor{light-gray}{gray}{0.70}

\title{Variable metric inexact line--search based methods for nonsmooth optimization
			 \thanks{This work has been partially supported by MIUR under the two projects FIRB - Futuro in Ricerca 2012, contract RBFR12M3AC and PRIN 2012, contract 2012MTE38N. Ignace Loris is a Research Associate of the Fonds de la Recherche Scientifique - FNRS. The Italian GNCS - INdAM is also acknowledged.}}

\author{S. Bonettini\footnotemark[1]
\and I. Loris\footnotemark[2]
\and F. Porta\footnotemark[1]
\and M. Prato\footnotemark[3]}

\begin{document}

\maketitle
\renewcommand{\thefootnote}{\fnsymbol{footnote}}

\footnotetext[1]{Dipartimento di Matematica e Informatica, Universit\`a di Ferrara, Via Saragat 1, 44122 Ferrara, Italy ({\tt silvia.bonettini@unife.it,federica.porta@unife.it}).}
\footnotetext[2]{D\'epartement de Math\'ematique, Universit\'e Libre de Bruxelles, Boulevard du Triomphe, 1050 Bruxelles, Belgium ({\tt igloris@ulb.ac.be}).}
\footnotetext[3]{Dipartimento di Scienze Fisiche, Informatiche e Matematiche, Universit\`a di Modena e Reggio Emilia, Via Campi 213/b, 41125 Modena, Italy
        				({\tt marco.prato@unimore.it}).}

\renewcommand{\thefootnote}{\arabic{footnote}}

\begin{abstract}
We develop a new proximal--gradient method for minimizing the sum of a differentiable, possibly nonconvex, function plus a convex, possibly non differentiable, function. The key features of the proposed method are the definition of a suitable descent direction, based on the proximal operator associated to the convex part of the objective function, and an Armijo--like rule to determine the step size along this direction ensuring the sufficient decrease of the objective function. In this frame, we especially address the possibility of adopting a metric which may change at each iteration and an inexact computation of the proximal point defining the descent direction. For the more general nonconvex case, we prove that all limit points of the iterates sequence are stationary, while for convex objective functions we prove the convergence of the whole sequence to a minimizer, under the assumption that a minimizer exists. In the latter case, assuming also that the gradient of the smooth part of the objective function is Lipschitz, we also give a convergence rate estimate, showing the ${\mathcal O}(\frac 1 k)$ complexity with respect to the function values. We also discuss verifiable sufficient conditions for the inexact proximal point and we present the results of a numerical experience on a convex total variation based image restoration problem, showing that the proposed approach is competitive with another state-of-the-art method.
\end{abstract}

\begin{keywords}
Proximal algorithms, nonsmooth optimization, generalized projection, nonconvex optimization.
\end{keywords}

\begin{AMS}
65K05, 90C30
\end{AMS}

\pagestyle{myheadings}
\thispagestyle{plain}
\markboth{S.~BONETTINI ET AL.}{VARIABLE METRIC INEXACT LINE--SEARCH BASED METHODS}

\section{Introduction}

In this paper we consider the problem
\begin{equation}\label{minf}
\min_{\x\in\R^n} f(\x) \equiv f_0(\x) + f_1(\x)
\end{equation}
where $f_1$ is a proper, convex, lower semicontinuous function and $f_0$ is smooth, i.e. continuously differentiable, on an open subset $\Omega_0$ of $\R^n$ containing $\dom(f_1)=\{\x\in\R^n: f_1(\x) <+\infty\}$.

We also assume that $f_1$ is bounded from below and that $\dom(f_1)$ is non-empty and closed. 
Formulation \eqref{minf} includes also constrained problems over convex sets, which can be introduced by adding to $f_1$ the indicator function of the feasible set.

When in particular $f_1$ reduces to the indicator function of a convex set $\Omega$, i.e. $f_1 = \iota_\Omega$ with
\begin{equation*}
\iota_\Omega(\x) = \left\{\begin{array}{cl} 0&\mbox{ if } \x\in\Omega \\ +\infty & \mbox{ if } \x\not\in\Omega.\end{array}\right.,
\end{equation*}
a simple and well studied algorithm for the solution of \eqref{minf} is the gradient projection (GP) method, which is particularly appealing for large scale problems. In the last years, several variants of such method have been proposed \cite{Birgin-etal-2003,Bonettini-etal-2009,Duchi-Hazan-Singer-2011,Hager-etal-2009}, with the aim to accelerate the convergence which, for the basic implementation, can be very slow. In particular, reliable acceleration techniques have been proposed for the so called gradient projection method with line--search along the feasible direction \cite[Chapter 2]{Bertsekas-1999}, whose iteration consists in
\begin{equation}\label{intro0}
\xkk = \xk+\lamk(\yk-\xk),
\end{equation}
where $\yk$ is the Euclidean projection of the point $\xk-\nabla f_0(\xk)$ onto the feasible set $\Omega$ and $\lamk\in[0,1]$ is a steplength parameter ensuring the sufficient decrease of the objective function. Typically, $\lamk$ is determined by means of a backtracking loop until an Armijo-type inequality is satisfied. Variants of the basic scheme are obtained by introducing a further variable stepsize parameter $\ak$, which controls the step along the gradient, in combination with a variable choice of the underlying metric. In practice, the point $\yk$ can be defined as
\begin{equation}\label{intro1}
\yk = \arg\min_{\y\in\Omega}\nabla f_0(\xk) ^T(\y-\xk) + \frac{1}{2\ak}(\y-\xk)^TD_k(\y-\xk)
\end{equation}
where $\ak$ is a positive parameter and $D_k\in \R^{n\times n}$ is a symmetric positive definite matrix. The stepsizes $\ak$ and the matrices $D_k$ have to be considered as ``free'' parameters of the method and a clever choice of them can lead to significant improvements in the practical convergence behaviour \cite{Birgin-etal-2003,Bonettini-etal-2013a,Bonettini-etal-2009}.

In this paper we generalize the GP scheme \eqref{intro0}--\eqref{intro1}, by introducing the concept of descent direction for the case where $f_1$ is a general convex function and we propose a suitable variant of the Armijo rule for the nonsmooth problem \eqref{minf}.
In particular, we focus on the case when the descent direction has the form $\yk-\xk$, with
\begin{equation}\label{intro2}
\yk = \arg\min_{\y\in\R^n} \nabla f_0(\xk) ^T(\y-\xk) + \dsk(\y,\xk) + f_1(\y)-f_1(\xk),
\end{equation}
where $\dsk(\cdot,\cdot)$ plays the role of a distance function, depending on the parameter $\sk\in \R^q$. Clearly, \eqref{intro2} is a generalization of \eqref{intro1}, which is recovered when $f_1 = \iota_\Omega$, by setting $\ds(\y,\x) = \frac 1\alpha (\y-\x)^TD(\y-\x)$, with $\sigma = (\alpha,D)$.

Formally, the scheme \eqref{intro0}-\eqref{intro2} is a forward--backward (or proximal gradient) method \cite{Combettes-Wajs-2005,Combettes-Pesquet-2011} depending on the parameters $\lamk$, $\sk$.

In particular, we deeply investigate the variant of the scheme \eqref{intro0}--\eqref{intro2} where the minimization problem in \eqref{intro2} is solved inexactly and we devise two types of admissible approximations. We show that both approximation types can be practically computed when $f_1(\x) = g(Ax)$, where $A\in \R^{m\times n}$ and $g:\R^m\to\R$ is a proper, convex, lower semicontinuous function with an easy-to-compute resolvent operator. In this case, our scheme consists in a double loop method, where the inner loop is provided by an implementable stopping criterion. For general $f_0$, we are able to prove that any limit point of the sequence generated by our inexact scheme is stationary for problem \eqref{minf}. The proof of this fact is essentially based on the properties of the Armijo-type rule adopted for computing $\lamk$ and it does not require any Lipschitz property of the gradient of $f_0$. When $f_0$ is convex, we prove a stronger result, showing that the iterates converge to a minimizer of \eqref{minf}, if it exists. In the latter case, under the further assumption that $\nabla f_0$ is Lipschitz continuous, we give a ${\mathcal O}(\frac 1 k)$ convergence rate estimate for the objective function values. Our analysis includes as special cases several state-of-the-art methods, as those in \cite{Birgin-etal-2003,Bonettini-Prato-2015a,Bonettini-etal-2009,Loris-Porta-2015,Tseng-Yun-2009}.

Forward--backward algorithms based on a variable metric have been recently studied also in \cite{Combettes-Vu-2014} for the convex case and in \cite{Chouzenoux-etal-2014} for the nonconvex case under the Kurdyka-{\L}ojasiewicz assumption (see also \cite{Frankel-etal-2015}). Even if our scheme is formally very similar to those in \cite{Chouzenoux-etal-2014,Combettes-Vu-2014}, the involved parameters have a substantially different meaning. In our case, the theoretical convergence is ensured by the Armijo parameter $\lamk$ in combination with the descent direction properties; this results in an almost complete freedom to choose the other algorithm parameters (e.g. $\ak$ and $D_k$), without necessarily relating them to the Lipschitz constant of $\nabla f_0$ (actually, our analysis, except the convergence rate estimate, is performed without this assumption). We believe that this is also one of the main strength of our method, since acceleration techniques based on suitable choices of $\ak$ and $D_k$, originally proposed for smooth optimization, can be adopted, leading to an improvement of the practical performances. The other crucial ingredient of our method is the inexact computation of the minimizer in \eqref{intro2}: this issue has been considered in several papers in the context of proximal and proximal gradient methods (see for example \cite{Attouch-etal-2013,Chouzenoux-etal-2014,Salzo-Villa-2012,Villa-etal-2013} and references therein). The approach we follow in this paper is more similar to the one proposed in \cite{Villa-etal-2013} and has the advantage to provide an implementable condition for the approximate computation of the proximal point. Moreover, we also generalize the ideas proposed in \cite{Birgin-etal-2003} for the inexact computation of the projection onto a convex set. Finally, we also mention the papers \cite{Auslender-etal-2007,Auslender-Teboulle-2006,Auslender-Teboulle-2009,Eckstein-1993} for the use of non Euclidean distances in the context of forward--backward and proximal methods.

The paper is organized as follows: some background material is collected in Section \ref{sec:definitions}, while the concept of descent direction for problem \eqref{minf} is presented and developed in Section \ref{sec:descent}. In Section \ref{sec:algorithm}, the modified Armijo rule is discussed. Then, a general convergence result for line--search descent algorithms based on this rule is proved, in the nonconvex case. Two different inexactness criteria, called of $\epsilon$-type and $\eta$-type are proposed in Sections \ref{sec:epsilon-approx} and \ref{sec:eta-approx}, and the related implementation is discussed in Sections \ref{sec:eta-implementation} and \ref{sec:eps-implementation}. Section \ref{sec:convergence} deals with the convex case, where the convergence of an $\epsilon$-approximation based algorithm is proved and the related convergence rate is analyzed. The results of a numerical experience on a total variation based image restoration problem are presented in Section \ref{sec:numerical} while our conclusions are given in Section \ref{sec:conclusions}.

\paragraph{Notation} We denote the extended real numbers set as $\bR = \R\cup \{-\infty,+\infty\}$ and by $\R_{\geq0}$, $\R_{>0}$ the set of non-negative and positive real numbers, respectively. The scaled Euclidean norm of an $n$-vector $\x$, associated to a symmetric positive definite matrix $D$ is $\|\x\|_D= \sqrt{\x^TD\x}$. Given $\mu\geq 1$, we denote by ${\mathcal M}_\mu$ the set of all symmetric positive definite matrices with all eigenvalues contained in the interval $[\frac 1 \mu,\mu]$. For any $D\in {\mathcal M}_\mu$ we have that $D^{-1}$ also belongs to ${\mathcal M}_\mu$ and
\begin{equation}\label{ine_norm}
\frac 1\mu\|\x\|^2 \leq \|\x\|^2_{D}\leq \mu\|x\|^2
\end{equation}
for any $\x\in\R^n$.

\section{Definitions and basic properties}\label{sec:definitions}

We recall the following definitions.
\begin{Def} \cite[p.213]{Rokafellar-1970} Let $f$ be any function from $\R^n$ to $\bR$. The \emph{one sided directional derivative} of $f$ at $\x$ with respect to a vector $\d$ is defined as
\begin{equation}\label{dir-der}
f'(\x;\d) = \lim_{\lambda \downarrow 0} \frac{f(\x+\lambda \d)-f(\x)}{\lambda}
\end{equation}
if the limit on the right-hand side exists in $\bR$.
\end{Def}
When $f$ is smooth at $\x$, then $f'(\x;\d) = \nabla f(\x)^T\d$. When $f$ is convex, its directional derivative has the following property.
\begin{Thm}\label{teo:directional-convex}\cite[Theorem 23.1]{Rokafellar-1970}
If $f$ is convex and $\x\in\dom(f)$, then for any $\d\in\R^n$ the limit at the right-hand side of \eqref{dir-der} exists and $f'(\x;\d) = \inf_{\lambda > 0} \frac{f(\x+\lambda \d)-f(\x)}{\lambda}$.
\end{Thm}
As a consequence of the previous theorem, for any convex function $f$ we have that $f'(\x;\d)$ exists for any $\x\in\dom(f)$, $\d\in\R^n$ and
\begin{equation}\label{ine_dir_der}
f'(\x;\d) \leq f(\x+ \d)-f(\x).
\end{equation}
\begin{Def}\cite[p. 394]{Tseng-Yun-2009}
A point $\x$ is \emph{stationary} for problem \eqref{minf} if $\x\in\dom(f)$ and
\begin{equation}\label{stationary}
f'(\x;\d)\geq 0 \ \ \ \forall \d\in \R^n.
\end{equation}
\end{Def}
\begin{Def}\cite[\S 2.3]{Frankel-etal-2015}
The \emph{proximity} or \emph{resolvent} operator associated to a convex function $f:\R^n\rightarrow \bR$ in the metric induced by a symmetric positive definite matrix $D$ is defined as
$$\prox_{f}^D(\x) = \arg\min_{\z\in\R^n} f(\z) + \frac{1}{2}\|\z-\x\|^2_D,       \ \ \ \forall \x\in \R^n.$$
\end{Def}
We remark that $\prox_{f}^D$ is a Lipschitz continuous function whose Lipschitz constant is $\|D\|$.
\begin{Def}\label{def:dual}
Let $f:\R^n\rightarrow \bR$ be a convex function. The \emph{conjugate function} of $f$ is the function $f^*:\R^n\rightarrow \bR$ defined as $f^*(\y) = \sup_{\x\in\R^n} \x^T\y - f(\x)$ $\forall \y\in\R^n$.
\end{Def}
The following proposition states a useful property of the conjugate.
\begin{Pro}\label{pro:duality}
Let $f:\R^n\rightarrow \bR$, $g:\R^m\rightarrow \bR$ be two convex functions, $A\in \R^{m\times n}$.
If $f(x) = g(Ax)$, then $f^*(A^T\y)\leq g^*(\y)$ $\forall \y \in \R^m$. 
\end{Pro}
\proof By Definition \ref{def:dual} we have
\begin{equation*}
f^*(A^T\y)= \sup_{\x\in\R^n}\ \x^TA^T\y - f(\x)= \sup_{\x\in\R^n}\ (A\x)^T\y - g(A\x) = \!\!\sup_{\z\in\R^m,\z = A\x}\!\! \z^T\y-g(\z) \leq \sup_{\z\in\R^m} \z^T\y-g(\z) = g^*(\y) .
\end{equation*}
\endproof
\vspace{-8mm}
\begin{Def}\cite[p. 82]{Zalinescu-2002}
Given $\epsilon\in\R_{\geq 0}$, the $\epsilon$-subdifferential of a convex function $f:\R^n\rightarrow\bR$ at a point $\z\in\R^n$ is the set
\begin{equation}\label{eps-subdiff}
\partial_\epsilon f(\z) = \{\w\in \R^n: f(\x) \geq f(z) + (\x-\z)^Tw -\epsilon, \ \ \forall x\in \R^n\}.
\end{equation}
\end{Def}
If $\z\in\dom(f)$, then $\partial_\epsilon f(\z)\neq \emptyset$. For $\epsilon = 0$ the usual subdifferential set $\partial f(\z)$ is recovered.
A useful property of the $\epsilon$-subdifferential is the following one.
\begin{Pro}\label{pro:zalinescu}\cite[Theorem 2.4.4 (iv)]{Zalinescu-2002}
Let $f:\R^n\rightarrow\bR$ be a convex, proper, lower semicontinuous function. Then for any $\epsilon\in\R_{\geq 0}$ and for any $\x\in\R^n$ we have $\x^*\in \partial_\epsilon f(\x) \Leftrightarrow \x\in\partial_\epsilon f^*(\x^*)$.
\end{Pro}

\section{A family of descent directions}\label{sec:descent}

When $f$ is smooth, a vector $\d\in\R^n$ is said a descent direction for $f$ at $\x$ when $\nabla f(\x)^T d < 0$. In the nonsmooth case \eqref{minf}, we give the following definition, based on the directional derivative.
\begin{Def}\label{def:descent-direction}A vector $\d\in\R^n$ is a \emph{descent direction} for $f$ at $\x\in\dom(f)$ if $f'(\x;\d)< 0$.
\end{Def}
Thanks to Theorem \ref{teo:directional-convex}, the previous definition is well posed. In this section we define a family of descent directions for problem \eqref{minf}. To this end, we define the following set of non--negative functions.

Given a convex set $\Omega\subseteq \R^n$ and a set of parameters $S\subseteq \R^q$, we denote by ${\mathcal D}(\Omega,S)$ the set of any \emph{distance--like} function $d_\s:\R^n\times\R^n \rightarrow \R_{\geq 0} \cup\{+\infty\}$ continuously depending on $\s\in S$ such that for all $\z,\x\in\Omega$ we have:
\begin{enumerate}[label=$({\mathcal D}_{\arabic*})$]
\item \label{distcontinuous} $\ds(\z,\x)$ is continuous in $(\sigma,\z,\x)$;
\item \label{distsmooth}$\ds(\z,\x)$ is smooth w.r.t. $\z\in\Omega$;
\item \label{diststrongcvx}$\ds(\z,\x)$ is strongly convex w.r.t. $\z$:
\begin{equation*}
\ds(\z_2,\x)\geq \ds(\z_1,\x)+\nabla_1\ds(\z_1,\x)^T(\z_2-\z_1)+\frac{m}{2}\|\z_2-\z_1\|^2\qquad\forall \z_1,\z_2\in\Omega,
\end{equation*}
where $m>0$ does not depend on $\sigma$ or $\x$ (here $\nabla_1$ denotes the gradient with respect to the first argument of a function);
\item \label{distzero} $\ds(\z,\x)= 0$  if and only if $\z=\x$ (which implies that $\nabla_1 \ds(\x,\x) = 0$ for all $\x\in \Omega$).
\end{enumerate}

The scaled Euclidean distance
\begin{equation}\label{scaled-euclidean}
\ds(\x,\y) = \frac{1}{2\alpha}\|\x-\y\|^2_{D}
\end{equation}
with $\vesi = (\alpha,D)$, where $\alpha > 0$ and $D\in\R^{n\times n}$ is a symmetric positive definite matrix, is an interesting example of a function in ${\mathcal D}(\R^n,S)$. Other examples of distance--like functions can be obtained by considering Bregman distances associated to a strongly convex function.

For a given array of parameters $\s\in S\subseteq\R^q$, let us introduce the function $\hs:\R^n\times\R^n \rightarrow \bR$ defined as
\begin{equation}\label{hxy}
\hsxy\z\x = \nabla f_0(\x)^T(\z-\x) + d_\sigma(\z,\x) + f_1(\z)-f_1(\x)\ \ \forall \z,\x\in\R^n,
\end{equation}
where $d_\s\in {\mathcal D}(\Omega,S)$ and $\Omega = \dom(f_1)$. We remark that $\hs$ depends continuously on $\s$, as $\ds$ does. Moreover, since $\ds(\cdot,\x)$ and $f_1$ are convex, proper and lower semicontinuous, $\hs(\cdot,\x)$ is also convex, proper and lower semicontinuous for all $\x\in\Omega_0$. Finally, for any point $\x\in\Omega$ and for any $\d\in \R^n$ we have
\begin{equation}\label{H2}
\dhsxy\x\d = f'(\x;\d),
\end{equation}
where $\dhs(\z,\x;\d)$ denotes the directional derivative of $\hsxy{\puntino}\x$ at the point $\z$ with respect to $d$. 
From assumption \ref{diststrongcvx}, it follows that $\hsxy{\cdot}{\x}$ is strongly convex and admits a unique minimum point for any $\x\in \Omega $.

Now we introduce the following operator $\p:\Omega_0\rightarrow\Omega$ associated to any function $\hs$ of the form \eqref{hxy}
\begin{equation}\label{proj}
\p(\x;\hs) = \arg\min_{\z\in\R^n} \hsxy \z\x.
\end{equation}
When $\ds$ is chosen as in \eqref{scaled-euclidean}, the operator \eqref{proj} becomes
\begin{equation*}
p(\x;\hs) = \prox_{\alpha f_1}^D(\x-\alpha D^{-1}\nabla f_0(\x)).
\end{equation*}
Under assumption \ref{diststrongcvx}, one can show that $\p(\x;\hs)$ depends continuously on $(\x,\sigma)$.
\begin{Pro}
Let $\ds\in {\mathcal D}(\Omega,S)$ and $\hs$ be defined as in \eqref{hxy}. Then $\p(\x;\hs)$ depends continuously on $(\x,\sigma)$.
\end{Pro}
\proof
Let $\y=\arg\min_{\z\in\R^n} \hsxy \z\x$. Then $\y$ is characterized by the equation $\nabla f_0(x)+\nabla_1 \ds(\y,\x) +w=0$, where $w\in\partial f_1(y)$. It follows that $f_1(u)\geq f_1(y)+w^T(u-y)$ for all $u\in\R^n$ or:
\begin{equation*}
f_1(u)\geq f_1(y)-(\nabla f_0(x)+\nabla_1 \ds(\y,\x))^T(u-y)\qquad\forall u\in \R^n.
\end{equation*}
Assumption \ref{diststrongcvx} expressed in $\y$ and $\u$ gives:
\begin{equation*}
\ds(u,\x)\geq \ds(\y,\x)+\nabla_1\ds(\y,\x)^T(u-\y)+\frac{m}{2}\|\y-u\|^2\qquad\forall u\in\R^n.
\end{equation*}
Together, these two inequalities yield:
\begin{equation*}
\frac{m}{2}\|\y-u\|^2\leq f_1(u)-f_1(\y)+\ds(u,x)-\ds(\y,\x)+\nabla f_0(x)^T(u-y)\qquad\forall u\in\R^n.
\end{equation*}
Let $\y_1=\p(\x_1;h_{\sigma_1})$ and $\y_2=\p(\x_2;h_{\sigma_2})$. Adding the previous inequality for $\y=\y_1$ (resp. $\y=\y_2$) and choosing $u=\y_2$ (resp. $u=\y_1$), one finds:
\begin{equation*}
m\|\y_1-\y_2\|^2\leq d_{\sigma_1}(\y_2,x_1)-d_{\sigma_1}(\y_1,\x_1)
+d_{\sigma_2}(\y_1,\x_2)-d_{\sigma_2}(\y_2,\x_2)+(\nabla f_0(\x_1)-\nabla f_0(\x_2))^T(\y_2-\y_1)
\end{equation*}
and hence:
\begin{equation*}
m\|\y_1-\y_2\|^2\leq d_{\sigma_2}(\y_1,\x_2)-d_{\sigma_1}(\y_1,\x_1)
+d_{\sigma_1}(\y_2,x_1)-d_{\sigma_2}(\y_2,\x_2)+\|\nabla f_0(\x_1)-\nabla f_0(\x_2)\|\,\|\y_2-\y_1\|.
\end{equation*}
It follows that $0\leq\|\y_1-\y_2\|\leq (b+\sqrt{b^2+4cm})/2m$ where $b=\|\nabla f_0(\x_1)-\nabla f_0(\x_2)\|$ and $c=d_{\sigma_2}(\y_1,\x_2)-d_{\sigma_1}(\y_1,\x_1)
+d_{\sigma_1}(\y_2,x_1)-d_{\sigma_2}(\y_2,\x_2)$. As $f_0$ is $C^1$, one has $\lim_{\x_2\to\x_1}b=0$. As $\ds(\z,\x)$ is continuous in $(\sigma,\z,\x)$, one also has that $\lim_{\x_2\to\x_1}c=0$. This shows then that $\lim_{\x_2\to\x_1}\|\y_2-\y_1\|=0$, in other words $\p(\x_1;h_{\sigma_1})$ is continuous in $(\sigma_1,\x_1)$.\endproof
Given a function $\ds\in {\mathcal D}(\Omega,S)$, we introduce also the function $\ths : \R^n\times\R^n\rightarrow\bR$ defined as
\begin{equation}\label{thxy}
\thsxy\z\x = \nabla f_0(\x)^T(\z-\x) + \gamma d_\sigma(\z,\x) + f_1(\z)-f_1(\x)\ \ \forall \z,\x\in\R^n
\end{equation}
for some $\gamma \in [0,1]$. We have
\begin{equation}\label{IL1}
\thsxy{\y}{\x}\leq \hsxy{\y}{\x}\ \ \ \forall\x,\y\in\R^n
\end{equation}
and $\ths = \hs$ when $\gamma = 1$. In the following we will show that
\begin{itemize}
\item the stationarity condition \eqref{stationary} can be reformulated in terms of fixed points of the operator $\p(\puntino;\hs)$;
\item the negative sign of $\ths$ detects a descent direction.
\end{itemize}
To this purpose, we collect in the following proposition some properties of the function $\hs$ and the associated operator $\p(\puntino;\hs)$.
\begin{Pro}\label{Proinexact1}
Let $\sigma\in S\subseteq\R^q$, $\gamma\in[0,1]$, $\ds\in{\mathcal D}$ and $\hs$, $\ths$ be defined as in \eqref{hxy}, \eqref{thxy}, where $\ds \in \mathcal{D}(\Omega,S)$. If $\x\in\Omega$ and $\y = \p(\x;\hs)$, then:
\begin{itemize}
\item[(a)] $\thsxy{\x}{\x}= 0$;
\item[(b)] if $\z\in\R^n$ and $\thsxy{\z}{\x} < 0$, then $f'(\x;\z-\x) < 0$;
\item[(c)] $\thsxy{\y}{\x}\leq 0$ and $\thsxy{\y}{\x} = 0$ if and only if $\y = \x$;
\item[(d)] $f'(\x;\y-\x) \leq 0$ and the equality holds if and only if $\thsxy\y\x=0$ (if and only if $\x=\y$).
\end{itemize}
\end{Pro}\noindent
{\it Proof. } (a) is a direct consequence of definition \eqref{thxy} and condition $({\mathcal D}_3)$ on $\ds$.

(b) If $\thsxy{\z}{\x} < 0$, we have
\begin{equation*}
0 \geq -  \gamma d_\sigma(\z,\x) > \nabla f_0(\x)^T(\z-\x)  + f_1(\z)-f_1(\x) \geq \nabla f_0(\x)^T(\z-\x)  + f_1'(\x;\z-\x) = f'(\x;\z-\x),
\end{equation*}
where the second inequality follows from definition \eqref{thxy} of $\ths$ and the third one from  \eqref{ine_dir_der}.

(c) Since $\y$ is the minimum point of $\hsxy{\puntino}{\x}$, part (a) with $\gamma = 1$ yields $\hsxy{\y}{\x}\leq 0$ which, in view of \eqref{IL1}, gives $\thsxy{\y}{\x}\leq 0$. If $\y = \x$, part (a) implies $\thsxy{\y}{\x} = 0$. Conversely, assume $\thsxy{\y}{\x} = 0$. From inequality \eqref{IL1} we have $\hsxy\y\x\geq 0$. On the other side, since $y$ is the minimum point of $\hs(\cdot,\x)$, part (a) with $\gamma = 1$ implies $\hsxy\y\x\leq 0$. Thus $\hsxy\y\x= 0$ and since $\y$ is the unique minimizer of $\hsxy{\cdot}{\x}$, we can conclude that $\x = \y$.

(d) From (c) we have $\thsxy\y\x\leq 0$. When $\thsxy\y\x<0$ then part (b) implies $f'(\x;\y-\x)<0$. When $\thsxy\y\x=0$, from (c) we obtain $\y=\x$ and, therefore, $f'(\x;\y-\x)=0$. Conversely, assume $f'(\x;\y-\x)=0$. This implies
\begin{equation*}
0 = \nabla f_0(\x)^T(\y-\x)  + f_1'(\x;\y-\x) \leq \nabla f_0(\x)^T(\y-\x)  + f_1(\y)-f_1(\x) \leq \thsxy\y\x.
\end{equation*}
Since $\thsxy\y\x\leq 0$, we necessarily have  $\thsxy\y\x=0$.
\endproof
The following proposition completely characterizes the stationary points of \eqref{minf} in two equivalent ways, as fixed points of the operator $\p(\cdot;\hs)$, i.e. the solutions of the equation $\x = \p(\x;\hs)$, or as roots of the composite function $r_{\vesi,\gamma}(\x) = \ths(p(\x;\hs),\x)$.
\begin{Pro}\label{Proproj}
Let $S \subseteq \R^q$, $\vesi\in S$, $\hs$, $\ths$ be defined as in \eqref{hxy}, $\gamma \in [0,1]$, $\ve x\in \Omega$ and $\y=\p(\x;\hs)$. The following statements are equivalent:
\begin{itemize}
\item[(a)] $\x$ is stationary for problem \eqref{minf};
\item[(b)] $\x=\y$;
\item[(c)] $\thsxy\y\x = 0$.
\end{itemize}
\end{Pro}\noindent
{\it Proof. }(a) $\Longleftrightarrow$ (b) Assume that $\x=\y$.
Then, $\hs(\cdot,\x)$ achieves its minimum at $x$ and inequality \eqref{stationary} applied to it yields $\dhsxy{\x}{\z-\x} \geq 0 \quad \forall \z\in \R^n$. Recalling \eqref{H2} we have $\dhsxy{\x}{\z-\x} = f'(\x;\z-\x)$, hence $\x$ is a stationary point for problem \eqref{minf}.

Conversely, let $\x\in\Omega$ be a stationary point of \eqref{minf} and assume by contradiction that $\x\neq\y$. Then, by Proposition \ref{Proinexact1} (d) we obtain $f'(\x,\y-\x)<0$, which contradicts the stationarity assumption on $\x$.

(b) $\Longleftrightarrow$ (c) See Proposition \ref{Proinexact1} (c).\hfill $\square$

\section{A line--search algorithm based on a modified Armijo rule}\label{sec:algorithm}

In this section we consider the modified Armijo rule described in Algorithm \ref{Algo2}, which is a generalization of the one in \cite{Tseng-Yun-2009}. Indeed the rule proposed in \cite{Tseng-Yun-2009} is recovered when $\ds$ is chosen as in \eqref{scaled-euclidean} and $\gamma\in[0,1)$. In the following we will prove that Algorithm \ref{Algo2} is well defined and classical properties of the Armijo condition still hold for this modified case.

\begin{algorithm}[ht]\caption{Modified Armijo linesearch algorithm}\label{Algo2}
Let $\{\xk\}_{k \in \N}$, $\{\tyk\}_{k \in \N}$ be two sequences of points in $\Omega$, and $\{\sk\}_{\kinN}$ be a sequence of parameters in $S$.   
Choose some $\delta, \beta \in(0,1)$, $\gamma \in [0,1]$. For all $k\in\N$ compute $\lambda^{(k)}$ as follows:
\begin{itemize}
\item[1.] Set $\lambda^{(k)}=1$ and $\dk=\tyk-\xk$.
\item[2.] \textsc{If}
\vspace{-5mm}
\begin{equation}\label{Armijo}
f(\xk+\lambda^{(k)}\dk)\leq f(\xk)+\beta \lambda^{(k)}\Delta\k
\end{equation}
\vspace{-3mm}
where
\vspace{-1mm}
\begin{equation}\label{Deltak}
\Delta\k = \thskxy\tyk\xk
\end{equation}
\vspace{-3mm}
\hspace*{-2.5mm}
\textsc{Then} go to step 3.\\\\
\textsc{Else} set $\lambda^{(k)}=\delta\lambda^{(k)}$ and go to step 2.
\item[3.] \textsc{End}
\end{itemize}
\end{algorithm}

Here and in the following we will define the function $\hsxy\cdot\cdot$ as in \eqref{hxy} and, for sake of simplicity, we will make the following assumption
\begin{itemize}
\item[(H0)] $\ds\in D(\Omega,S)$, where $\Omega=\dom(f_1)$ and $S\subseteq\R^q$ is a compact set.
\end{itemize}
\begin{Pro}\label{ProArm}
Let $\{\xk\}_\kinN$, $\{\tyk\}_\kinN$ be two sequences of points in $\Omega$, $\{\sk\}_\kinN$ a sequence of parameters in $S\subseteq \R^q$ and $\gamma\in[0,1]$.
Assume that
\begin{equation}\label{A2}
\thskxy\tyk\xk < 0 
\end{equation}
for all $k$. Then, the line--search Algorithm \ref{Algo2} is well defined, i.e. for each $k\in\N$ the loop at step 2 terminates in a finite number of steps. If, in addition, we assume that $\{\xk\}_{\kinN}$ and $\{\tyk\}_\kinN$ are bounded sequences and
$f(\xkk) \leq f(\xk)$, then we have that $\Delta\k= \thskxy\tyk\xk$ is bounded.
Assuming also that
\begin{equation}\label{A3}
\displaystyle\lim_{k\rightarrow\infty}f(\xk)- f(\xk+\lambda^{(k)}\ve d^{(k)})=0,
\end{equation}
where $\lambda^{(k)}$ and $\dk$ are computed with Algorithm \ref{Algo2}, then we have
\begin{equation*}
\lim_{k\rightarrow \infty} \thskxy\tyk\xk = 0.
\end{equation*}
\end{Pro}
{\it Proof.} 
We prove first that the loop at step 2 of Algorithm \ref{Algo2} terminates in a finite number of steps for any $k\in \N$.
Assume by contradiction that there exists a $k\in\N$ such that Algorithm \ref{Algo2} performs
an infinite number of reductions, thus, for any $j\in\N$, we have
\begin{eqnarray*}
\beta\Delta\k&<&\frac{f(\ve x^{(k)}+\delta^j\ve d^{(k)})-f(\ve x^{(k)})}{\delta^j}\\
&=&\frac{f_0(\ve x^{(k)}+\delta^j\ve d^{(k)})-f_0(\ve x^{(k)})}{\delta^j} + \frac{f_1(\ve x^{(k)}+\delta^j\ve d^{(k)})-f_1(\ve x^{(k)})}{\delta^j}\\
&\leq& \frac{f_0(\ve x^{(k)}+\delta^j\ve d^{(k)})-f_0(\ve x^{(k)})}{\delta^j} + \frac{\delta^jf_1(\ve x^{(k)}+\ve d^{(k)})+(1-\delta^j)f_1(\xk)-f_1(\ve x^{(k)})}{\delta^j}\\
&=&\frac{f_0(\ve x^{(k)}+\delta^j\ve d^{(k)})-f_0(\ve x^{(k)})}{\delta^j} +f_1(\tyk)-f_1(\xk),
\end{eqnarray*}
where the second inequality is obtained by means of the Jensen inequality applied to the convex function $f_1$.
Taking limits on the right hand side for $j\rightarrow\infty$ we obtain
\begin{eqnarray*}
\beta\Delta\k&\leq& \nabla f_0(\xk)^T\dk + f_1(\tyk)-f_1(\xk)\\
                 &\leq & \nabla f_0(\xk)^T\dk + f_1(\tyk)-f_1(\xk) + \gamma\dsk(\tyk,\xk)\\
                 &=   & \Delta\k  <   0,
\end{eqnarray*}
where the second inequality follows from the non--negativity of $\ds\in{\mathcal D}(\Omega,S)$ and the last one from \eqref{A2}. Since $0<\beta<1$, this is an absurdum.

Assume now that $\{\xk\}_\kinN$, $\{\tyk\}_\kinN$ are bounded sequences and that $f(\xkk)\leq f(\xk)$. We show that $\Delta\k=\thskxy\tyk\xk$ is bounded.
By assumption \eqref{A2}, $\thskxy\tyk\xk$ is bounded from above. We show that it is also bounded from below. Indeed we have
\begin{eqnarray*}
\thskxy\tyk\xk & = & \nabla f_0(\xk)^T(\tyk-\xk) +\gamma d_{\sigma^{(k)}}(\tyk,\xk) + f_1(\tyk)-f_1(\xk)\\
 & {\geq} & \nabla f_0(\xk)^T(\tyk-\xk) + f_1(\tyk)-f_1(\xk)\\
  & {=} & \nabla f_0(\xk)^T(\tyk-\xk) + f_1(\tyk)-f(\xk)+f_0(\xk)\\
  & {\geq} & \nabla f_0(\xk)^T(\tyk-\xk) + f_1(\tyk)-f(\x^{(0)})+f_0(\xk),
\end{eqnarray*}
where the first inequality follows from the non--negativity of $\ds$, the second one is obtained by adding and subtracting $f_0(\xk)$ and the last one is a consequence of $f(\xkk)\leq f(\xk)$.

As $f_1$ is proper and convex, there exists a supporting hyperplane, i.e. $\exists a,b\in\R^n$ such that $f_1(u)\geq a^T u+b$ for all $u\in\R^n$. Thus:
\begin{equation*}
\thskxy\tyk\xk\geq \nabla f_0(\xk)^T(\tyk-\xk)+ a^T \tyk+b-f(\x^{(0)})+f_0(\xk).
\end{equation*}
The right hand side is a continuous function of $\xk$ and $\tyk$. As these are assumed to lie on a closed and bounded set, the left hand side is bounded (from below) as well.

Let us show that the only limit point of $\Delta\k$ is zero. We observe that from \eqref{A2} and \eqref{A3} we obtain
\begin{equation}\label{add9}
0=\displaystyle\lim_{k\to\infty}f(\xk)- f(\xk+
\lambda^{(k)}\ve d^{(k)}) = \beta \lim_{k\to\infty}\Delta\k\lamk.
\end{equation}
Assume that there exists a subset of indices $K\subseteq\N$ such that $\lim_{k\in K,k\to\infty}\Delta\k=\bar{\Delta}\in\R$, with $\bar\Delta < 0$. By \eqref{add9}, this implies that
\begin{equation}\label{bdd1}
\lim_{k\in K,k\to\infty}\lamk=0.\end{equation}
Denote by $\bar K\subseteq K$ a set of indices such that $\lim_{k\in \bar K, k\to\infty}\vesik = \bar\vesi$, $\lim_{k\in \bar K,k\to\infty}\xk=\bar\x$ and $\lim_{k\in \bar K,k\to\infty}\tyk=\tilde\y$ for some $\bar\vesi\in S$, $\bar\x,\tilde\y\in\Omega$.
From \eqref{bdd1} we have that for any sufficiently large index $k\in\bar K$, Algorithm \ref{Algo2} makes at least a reduction: this means that
$$\beta
({\lambda^{(k)}}/{\delta}) \Delta\k<f(\ve x^{(k)}+({\lambda^{(k)}}/{\delta})\ve d^{(k)})-f(\ve x^{(k)}),$$
for all sufficiently large $k\in \bar K$.
Repeating the same arguments employed in the first part of the proof, we obtain
\begin{eqnarray*}
\beta
 \Delta\k&<&\frac{f_0(\ve x^{(k)}+({\lambda^{(k)}}/{\delta})\ve d^{(k)})-f_0(\ve x^{(k)})}{{\lambda^{(k)}}/{\delta}}+
 f_1(\tyk)-f_1(\xk)\\
 &\leq&\frac{f_0(\ve x^{(k)}+({\lambda^{(k)}}/{\delta})\ve d^{(k)})-f_0(\ve x^{(k)})}{{\lambda^{(k)}}/{\delta}}+
 f_1(\tyk)-f_1(\xk)+\gamma\ds(\tyk,\xk).
\end{eqnarray*}
Taking limits on both sides for $k\in \bar K,k\to\infty$, since $\{\ve d^{(k)}=\tyk-\xk\}_\kinN$ is bounded and
by \eqref{bdd1} we obtain $ \beta\bar{\Delta}\leq \bar{\Delta} <0  $,
which is an absurdum, being $0<\beta <1$.\endproof
We prove also the following useful Lemma.
\begin{Lemma}\label{Lem:1}
Let $\{\xk\}_\kinN$, $\{\tyk\}_\kinN$ be two sequences of points in $\Omega$, $\{\sk\}_\kinN$ a sequence of parameters in $S\subseteq \R^q$ and $\gamma\in[0,1]$.
Assume that
\begin{equation}\label{suff-decr}
f(\xkk) \leq f(\xk+\lamk\dk), \quad \dk=\tyk-\xk
\end{equation}
where $\tyk$ satisfies \eqref{A2} and $\lamk$ is computed by Algorithm \ref{Algo2} for any $k\in\N$. Suppose that $f$ is bounded from below. Then, we have
\begin{equation}\label{Lemma:1}
0\leq-\sum_{k=0}^\infty \lamk\thskxy\tyk\xk  <\infty.
\end{equation}
\end{Lemma}
\proof Denote by $\ell\in\R$ a lower bound for $f$, i.e. $\ell \leq f(\x)$ $\forall \x\in \R^n$. Inequalities \eqref{Armijo} and \eqref{suff-decr} can be combined as
\begin{equation*}
-\beta \lamk \thskxy\tyk\xk \leq f(\xk)-f(\xkk).
\end{equation*}
Summing the previous inequality for $k=0,...,j$ gives
\begin{equation}
-\beta \sum_{k=0}^j\lamk\thskxy\tyk\xk \leq \sum_{k=0}^j (f(\xk)-f(\xkk))= f(\ve x^{(0)})-f(\xjj) \leq f(\ve x^{(0)})-\ell.\label{add4}
\end{equation}
Thus, inequality \eqref{Lemma:1} follows.
\endproof

\subsection{A class of line--search based algorithms}

Proposition \ref{ProArm} allows the convergence analysis of a wide class of descent methods based on the Armijo condition \eqref{Armijo}. The crucial ingredients of these methods are
\begin{itemize}
\item a descent direction $\dk = \tyk-\xk$, where $\tyk$ is a suitable approximation of the point $\p(\xk;\hs)$;
\item the sufficient decrease of the objective function between two successive iterations, which has to amount at least to $\lamk\thsxy\tyk\xk$, where $\lamk$ is determined by the backtracking procedure given in Algorithm \ref{Algo2}.
\end{itemize}
\begin{Thm}\label{teo-suff-decr}
Let $\{\xk\}_\kinN$, $\{\tyk\}_\kinN$ be two sequences of points in $\Omega$, $\{\sk\}_\kinN\subset S$ and $\gamma\in[0,1]$. Assume that there exists a limit point $\bar{\ve x}$ of $\{\xk\}_{k\in\N}$ and let $K'\subseteq \N$ be a subset of indices such that $\lim_{k\in K',k\rightarrow\infty} \xk = \barx\in\Omega$.
Assume that, for any $k\in\N$ we have
\begin{equation*}
f(\xkk)\leq f(\xk+\lamk\dk), \ \ \dk=\tyk-\xk,
\end{equation*}
where $\lamk$ is computed by Algorithm \ref{Algo2}, $\tyk$ satisfies \eqref{A2} and there exists $K''\subseteq K'$ such that
\begin{equation}\label{dist}\lim_{k\in K'',k\rightarrow\infty}\hskxy\tyk\xk-\hskxy{\yk}\xk=0, \ \ \mbox{ with } \ \ \yk=\p(\xk;\hsk).\end{equation}
Then $\bar{\ve x}$ is a stationary point for problem \eqref{minf}.
\end{Thm}
\proof First, we notice that Algorithm \ref{Algo2} is well defined, since \eqref{A2} holds. 
We observe that, since $\hsk$ is strongly convex with modulus of convexity $m$ and $\yk$ is its minimum point, we have
\begin{equation}\label{hstrongly_convex}
\frac m 2 \|\z - \yk\|^2 \leq \hskxy\z\xk-\hskxy\yk\xk \ \ \ \forall \z\in \R^n.
\end{equation}
Setting $\z = \tyk$ in the previous inequality and using \eqref{dist} gives
\begin{equation}\label{dist2}
\lim_{k\in K'',k\rightarrow\infty}\|\tyk - \yk\|=0.
\end{equation}
By continuity of the operator $\p(\x;\hs)$, since $\{\xk\}_{k\in K'}$ is bounded, $\{\yk\}_{k\in K'}$ is bounded as well. Thus, \eqref{dist2} implies that $\{\tyk\}_{k\in K''}$ is also bounded and there exists a limit point $\bar \y$ of $\{\tyk\}_\kinN$.
We define $K\subseteq K''$ such that $\lim_{k\in K,k\rightarrow\infty}\tyk = \bar\y$ and $\lim_{k\in K, k\rightarrow \infty} \sk = \bar{\s}$. By continuity of the operator $\p(\x;\hs)$ with respect to all its arguments, \eqref{dist2} implies that $\bar\y= p(\bar x;\hsbar)$.

Consider now the sequence $\{f(\xk)\}_{k\in\N}$. From assumption \eqref{suff-decr} it follows that
\begin{equation} \label{ine3}
f(\xkk)\leq f(\xk+\lamk\dk)\leq f(\xk).
\end{equation}
Thus, the sequence $\{f(\xk)\}_{k\in\N}$ is monotone nonincreasing and, therefore, it converges to some $\bar f\in\bR$. Since $f$ is lower semicontinuous and $\bar{\ve x}$ is a limit point of $\{\xk\}_{k\in\N}$, we have
\begin{equation*}
\bar f = \lim_{k\rightarrow\infty} f(\xk) = \lim_{k\rightarrow\infty} f(\xkk) \geq f(\bar{\ve x}).
\end{equation*}
The previous inequality implies that $\bar f\in \R$ and this fact, together with inequality \eqref{ine3}, gives
\begin{equation*}
\lim_{k\rightarrow\infty} f(\xk) - f(\xk+\lamk\dk) = 0 .
\end{equation*}
Thus we can apply Proposition \ref{ProArm} and obtain
\begin{equation*}
\lim_{k\rightarrow\infty,k\in K} \thskxy\tyk\xk = 0.
\end{equation*}
Combining the previous equality with \eqref{IL1} and \eqref{dist} yields
\begin{equation*}
0=\lim_{k\rightarrow\infty,k\in K} \thskxy\tyk\xk\leq \lim_{k\rightarrow\infty,k\in K} \hskxy\tyk\xk =\lim_{k\rightarrow\infty,k\in K} \hskxy\yk\xk.
\end{equation*}
Since $\hskxy\yk\xk\leq 0$, this implies $\lim_{k\rightarrow\infty,k\in K} \hskxy\yk\xk=0$. Expressing inequality \eqref{hstrongly_convex} for $\z = \xk$, we can write
\begin{equation}\nonumber
\frac m 2 \|\xk-\yk\|^2 \leq \hskxy\xk\xk-\hskxy\yk\xk = -\hskxy\yk\xk \stackrel{k\rightarrow\infty,k\in K}{\longrightarrow} 0.
\end{equation}
Thus, we proved that $\bar \y = \bar\x$,
and, by Proposition \ref{Proproj} we have that $\bar\x$ is stationary.
\endproof
Let us now discuss assumption \eqref{dist} in the previous theorem, concerning the inexact solution of the minimum problem in \eqref{proj}. Assumption \eqref{A2} guarantees that $\dk = \tyk-\xk$ is a descent direction, which is needed for the line--search algorithm. However, it is not sufficient to ensure that the limit points are stationary, but we need also to assume that \eqref{dist} holds.

As counterexample, consider the case $n=1$, $f_0(x)=x^2/2$, $f_1(x)=0$, $d_\sigma(x,y)=(x-y)^2/2$, $\beta=\delta=1/2$. The sequence $\xkk = \xk+\lamk(\tyk-\xk)$ with $\lamk=1$, $\tyk = \xk-(1/2)^{k+1}$ satisfies all the assumptions of Theorem \ref{teo-suff-decr} except \eqref{dist}. However, starting from $\x^{(0)}=2$, the sequence writes as $\xk = 1+(1/2)^k\stackrel{k\to\infty}{\rightarrow} 1$, while the only stationary point is 0.

We remark that assumption \eqref{dist} could be replaced by requiring that $f_1$ is continuous and \eqref{dist2} holds. Clearly, \eqref{dist} cannot be checked directly, but it is very general. In the following sections, we will consider two implementable conditions which imply \eqref{dist} and in Sections \ref{sec:eta-implementation}--\ref{sec:eps-implementation} we show how $\tyk$ can be computed in practice without knowing $\p(\xk;\hsk)$.

\subsection{$\epsilon$- approximations}\label{sec:epsilon-approx}

In this section we will assume that $\ds$ has the form \eqref{scaled-euclidean} and, in this case, we will describe a sufficient condition for \eqref{dist}.

We observe that $\y=\p(\x;\hs)=\prox_{\alpha f_1}^D(\x-\alpha D^{-1}\nabla f_0(\x))$ if and only if $0\in \partial \hsxy\y\x$, that is
\begin{equation}\label{opt-eps}
\frac{1}{\alpha}D(\z-\y)\in\partial f_1(\y),
\end{equation}
where $\z = \x -\alpha D^{-1}\nabla f_0(\x)$. Borrowing the ideas in \cite{Salzo-Villa-2012, Villa-etal-2013}, we consider a relaxed version of \eqref{opt-eps} and we study the properties of any point $\tilde\y$ satisfying the following inclusion
\begin{equation}\label{eps-approx}
\frac 1 \alpha D(\z -\tilde y) \in \partial_\epsilon f_1(\tilde y),
\end{equation}
where $\epsilon \in\R_{\geq 0}$.
\begin{Lemma}\label{lem:eps-approx}
Let $\ds$ be defined as in \eqref{scaled-euclidean} and $\x\in\Omega$. Assume that $\y = p(\x;\hs)$ and that $\ty$ satisfies \eqref{eps-approx} for some $\epsilon \in \R_{\geq 0}$. Then $\ty\in \Omega$ and we have
\begin{itemize}
\item[(a)] $\hsxy\ty\x-\hsxy\y\x\leq \epsilon$;
\item[(b)] $ \|\ty - \y \|^2\leq \alpha\mu\epsilon$, for all $\mu\in\R_{>0}$ with $\frac 1 \mu \leq \lambda_{\min}(D)$, $\lambda_{\min}$ being the smallest eigenvalue of $D$.
\end{itemize}
\end{Lemma}
\proof
Since we have $\partial_\epsilon \hsxy\ty\x \supseteq\{\frac 1 \alpha D(\ty-\z) + \w: \w\in\partial_\epsilon f_1(\ty) \}$ (see \cite[Theorem 2.4.2 viii]{Zalinescu-2002}), inclusion \eqref{eps-approx} implies $0\in \partial_\epsilon \hs(\ty,\x)$ which, by definition \eqref{eps-subdiff} of $\epsilon$-subdifferential, is equivalent to
\begin{equation}\label{funval}
\hsxy\w\x \geq \hsxy{\tilde y}\x-\epsilon\ \ \ \forall w\in \R^n.
\end{equation}
We recall that $\hs(\puntino,\x)$ is strongly convex with modulus $m=2/(\alpha\mu)$ and $\y$ is its minimizer. This yields
\begin{equation}\nonumber
\frac{1}{\alpha\mu} \|\tilde y - \y\|^2 \leq \hsxy{\tilde y}\x-\hsxy\y\x\leq \epsilon,
\end{equation}
where the rightmost inequality follows from \eqref{funval} with $\w=\y$.
%
\endproof
The previous result combined with Theorem \ref{teo-suff-decr} directly implies the following Corollary.
\begin{Cor}\label{cor:eps-approx}
Let $0<\alpha_{\min}\leq \alpha_{\max}$, $\gamma\in[0,1]$, $\mu \geq 1$. Assume that $\{\ak\}_\kinN\subset[\alpha_{\min},\alpha_{\max}]$, $\{D_k\}_\kinN\subset {\mathcal M}_\mu$, $\{\epsilon_k\}_\kinN\subset \R_{\geq 0}$, $\lim_{k\rightarrow\infty} \epsilon_k = 0$. Let $\{\xk\}_\kinN$, $\{\tyk\}_\kinN$ be two sequences of points in $\Omega$ such that, for any $k\in\N$, \eqref{suff-decr} holds, where $\lamk$ is computed by Algorithm \ref{Algo2} and $\tyk$ satisfies \eqref{A2} and
\begin{equation}\label{eps-approxk}
\frac 1 {\ak} D_k(\zk - \tyk ) \in \partial_{\epsilon_k} f_1(\tyk),
\end{equation}
with $\zk =  \xk-\ak D_k^{-1} \nabla f_0(\xk)$.
Then, any limit point of the sequence $\{\xk\}_\kinN$ is stationary for problem \eqref{minf}.
\end{Cor}

\subsection{$\eta$-approximations}\label{sec:eta-approx}

A different approach to define a suitable approximation of the operator \eqref{proj} is based on the following definition.
\begin{equation}\label{etaset}
P_\eta(\x;\hs) = \{\ty\in \Omega : \hs(\ty,\x)\leq \eta \hs(\y,\x),\mbox{ where } \y = p(\x;\hs) \}
\end{equation}
for some $\eta\in (0,1]$. This idea of inexactness was introduced first in \cite{Birgin-etal-2003} to approximate the projection operator onto a convex set in the context of scaled gradient projection methods for smooth optimization. Clearly, if
\begin{equation}\label{eta-approx0}\ty\in P_\eta(\x;\hs),\end{equation}
then $\hsxy\ty\x\leq 0$ and $\hsxy\ty\x= 0$ if and only if $\hsxy\y\x= 0$ which implies $\ty = \y$.

The following Theorem establishes a convergence result under the condition $\tyk\in P_\eta(\xk;\hs)$.
\begin{Thm}\label{thm:eta-approx}
Let $\eta\in(0,1]$, $0\leq\gamma\leq  1$, $\{\sk\}_\kinN\subset S$ and $\{\xk\}_\kinN\subset \Omega$ satisfying \eqref{suff-decr}, where $\lamk$ is computed by Algorithm \ref{Algo2}, with
\begin{equation}\label{eta-approxk}
\tyk \in P_\eta(\xk;\hsk).
\end{equation}
Then, either for some $k$ the iterate $\xk$ is stationary for problem \eqref{minf}, or any limit point $\bar{\ve x}$ of $\{\xk\}_\kinN$ is stationary for problem \eqref{minf}.
\end{Thm}
\proof
We set $\yk = \p(\xk;\hsk)$ and we first observe that $\gamma\leq 1$ and \eqref{eta-approxk} imply
\begin{equation}\label{add7}
\thskxy\tyk\xk\leq \hskxy\tyk\xk\leq \eta \hskxy\yk\xk\leq 0.
\end{equation}
If at some iterate $k\in\N$ we have $\tilde h_{\sk, \gamma}(\tyk,\xk)=0$ and, as a consequence, $ \hskxy\yk\xk=0$, then, by Proposition \ref{Proproj}, $\xk$ is a stationary point for problem \eqref{minf}.

Otherwise $\tilde h_{\sk, \gamma}(\tyk,\xk)< 0$ for all $k\in\N$ and, thus, \eqref{A2} holds. Consider now a limit point $\bar x\in\Omega $ of $\{\xk\}_{\kinN}$ (if one exists) such that $\lim_{k\rightarrow\infty,k\in K'} \xk = \bar x$ for some set of indices $K'\subseteq \mathbb N$.\\
We first prove that $\{\tyk\}_{k\in K'}$ is bounded, using the strong convexity of $\hsk(\cdot,\xk)$. From \eqref{eta-approxk} we have
\begin{equation}\label{add8}
\hskxy\tyk\xk- \hskxy\yk\xk \leq (\eta-1) \hskxy\yk\xk.
\end{equation}
Since $\hsk(\cdot,\xk)$ is strongly convex with modulus of convexity $m$, and $\yk$ is the minimizer of $\hsk(\cdot,\xk)$, we can write
\begin{equation}\nonumber
\frac{ m}{2}\|\tyk-\yk\|^2 \leq \hsk(\tyk,\xk)- \hsk(\yk,\xk) \leq (\eta-1) \hsk(\yk,\xk).
\end{equation}
Since $\yk$ depends continuously on $\xk$, when $\{\xk\}_{k\in K'}$ is bounded, and all lie in a closed set, then $\{\yk\}_{k \in K'}$ is also bounded. Recalling Proposition \ref{ProArm}, we have that $\{\thskxy\tyk\xk\}_{k\in K'}$ is bounded from below; then, using inequalities \eqref{add7}, we can conclude that $\hsk(\yk,\xk)$ is also bounded from below for $k\in K'$ and, thus, $\{\tyk\}_{k\in K'}$ is bounded.
We define $K\subseteq K'$ as the set of indices such that $\lim_{k\in K, k\rightarrow +\infty}\vesik = \bar\vesi$, $\lim_{k\in K,k\rightarrow+\infty}\yk=\bar\y$
for some $\bar\vesi\in S$, $\bar\y\in\Omega$. Thanks to the continuity of the operator \eqref{proj}, the set $K$ is well defined,  since the sequences $\{\xk\}_{k\in K'}$, $\{\vesik\}_\kinN$ are bounded, and, moreover, we have $\bar\y = p(\bar\x;h_{\bar\sigma})$.
Reasoning as in the proof of Theorem \ref{teo-suff-decr}, the existence of a limit point guarantees that \eqref{A3} is satisfied. Then, by Proposition \ref{ProArm}, we obtain $\lim_{k\rightarrow\infty,k\in K} \thskxy\tyk\xk =0$. 
Combining this with \eqref{eta-approxk}, we also have
\begin{equation*}
0=\lim_{k\rightarrow\infty,k\in K}\thskxy\tyk\xk \leq\lim_{k\rightarrow\infty,k\in K}\hsk(\tyk,\xk) \leq \eta \lim_{k\rightarrow\infty,k\in K} \hsk(\yk,\xk)
\end{equation*}
which, since $\hsk(\yk,\xk)\leq 0$, implies
\begin{equation}\label{add_new1}
\lim_{k\rightarrow\infty,k\in K} \hsk(\yk,\xk) = 0
\end{equation}
Invoking again the strong convexity of $\hsk(\puntino,\xk)$, we obtain
\begin{equation}\nonumber
\frac{ m}{2}\|\xk-\yk\|^2 \leq \hsk(\xk,\xk)- \hsk(\yk,\xk) = - \hsk(\yk,\xk)
\end{equation}
with, together with \eqref{add_new1} gives $\lim_{k\to\infty,k\in K} \|\yk-\xk\|^2 = 0$. Thus, $\bar\y=\bar\x$ and
by Proposition \ref{Proproj}, we conclude that $\barx$ is stationary.
\endproof

\subsection{Remarks}

Different notions of inexactness have been proposed in the literature (see \cite{Salzo-Villa-2012,Villa-etal-2013} and references therein), especially in the context of proximal point methods, with the aim of approximating the resolvent operator, and some of them could be considered also in our framework. A synthetic description of possible inexactness notions and their relationships is given in Figure \ref{fig:inexact}.

\begin{figure}[ht]
\begin{center}
\resizebox{\textwidth}{!}{%
\framebox{\begin{tabular}{ccccccc}
$\frac1\alpha D(z-\tilde \y) \in\partial_\epsilon f_1({\tilde \y})$ \\
&\rotatebox{-45}{$\Rightarrow$}\\
& &$0\in\partial_\epsilon \hsxy{\tilde\y}\x$&$\Leftrightarrow$&$\hsxy\ty\x\leq \hsxy\y\x+\epsilon$&$\Rightarrow$&$\|\ty-\y\|^2\leq \kappa \epsilon$\\
&\rotatebox{45}{$\Rightarrow$}\\
$\mbox{dist}(0,\partial \hsxy\ty\x) \leq \epsilon$\\
(when $D = I$)
\end{tabular}}}
\caption{Connection of different inexactness notions, under the assumption \eqref{scaled-euclidean}. The proof of the implications are given in Lemma \ref{lem:eps-approx} and in \cite[Proposition 1]{Salzo-Villa-2012}.}\label{fig:inexact}
\end{center}
\end{figure}
It is difficult to insert the inexactness criterion \eqref{eta-approx0} in the scheme in Figure \ref{fig:inexact}, since the shape of $P_\eta$ in \eqref{eta-approx0} depends on $\x$, while the implications in Figure \ref{fig:inexact} are independent of $\x$.\\
In general, we observe that from inequality \eqref{add8} and by definition of $\epsilon$-subdifferential we have
\begin{equation}\nonumber
0\in \partial_{\epsilon_k}  \hsk(\tyk,\xk), \mbox{ with } \ \epsilon_k = (\eta-1) \hsk(\yk,\xk).
\end{equation}
We give a pictorial example of the sets of admissible approximations $\ty$ of the exact minimizer $\y$ defined by conditions \eqref{eta-approx0} and \eqref{eps-approx} in Figure \ref{fig:inexact-pictures}. This example refers to the case where $f_1(\x) = \iota_\Omega(\x)$ is the indicator function of a convex closed set $\Omega\subseteq \R^n$. Choosing the Euclidean metric, i.e. \eqref{scaled-euclidean} with $D=I$, $\alpha = 1$, as distance function, the operator $\p(\x;\hs)$ reduces to the Euclidean projection of the point $\z = \x-\nabla f_0(x)$ onto $\Omega$. Moreover, condition \eqref{eps-approx} becomes
\begin{equation}\label{esempio}
\ty\in\Omega \mbox{ and } (\w-\ty)^T(\z-\ty) \leq \epsilon. \ \ \ \forall \w\in\Omega.
\end{equation}
As well explained in \cite{Salzo-Villa-2012,Villa-etal-2013}, from a geometrical point of view, a point $\ty\in\Omega$ satisfies \eqref{esempio} if and only if $\Omega$ is contained in the negative half-space determined by the hyperplane of equation $ (\w-\ty)^T(\z-\ty)/\|\z-\ty\| = \epsilon/\|\z-\ty\|$, which is normal to $\z-\ty$ at a distance $\epsilon/\|\z-\ty\|$ from $\ty$.

On the other side, setting $\gamma = 1$ for simplicity, we have $\ths(\puntino,\x)=\hs(\puntino,\x) = \frac 1 2 \|\puntino-\z\|^2-\frac 1 2 \|\x-\z\|^2 + \iota_\Omega(\puntino)-\iota_\Omega(\x)$. Thus, the set $P_\eta(\x;\hs)$ is the intersection of the set $\Omega$ with the ball centered in $\z$ of radius $\sqrt{\eta\|\y-\z\|^2+(1-\eta)\|\x-\z\|^2} $.
\def\FigScala{0.4}
\begin{figure}
\begin{center}
\begin{tabular}{cc}
{\includegraphics[scale=\FigScala]{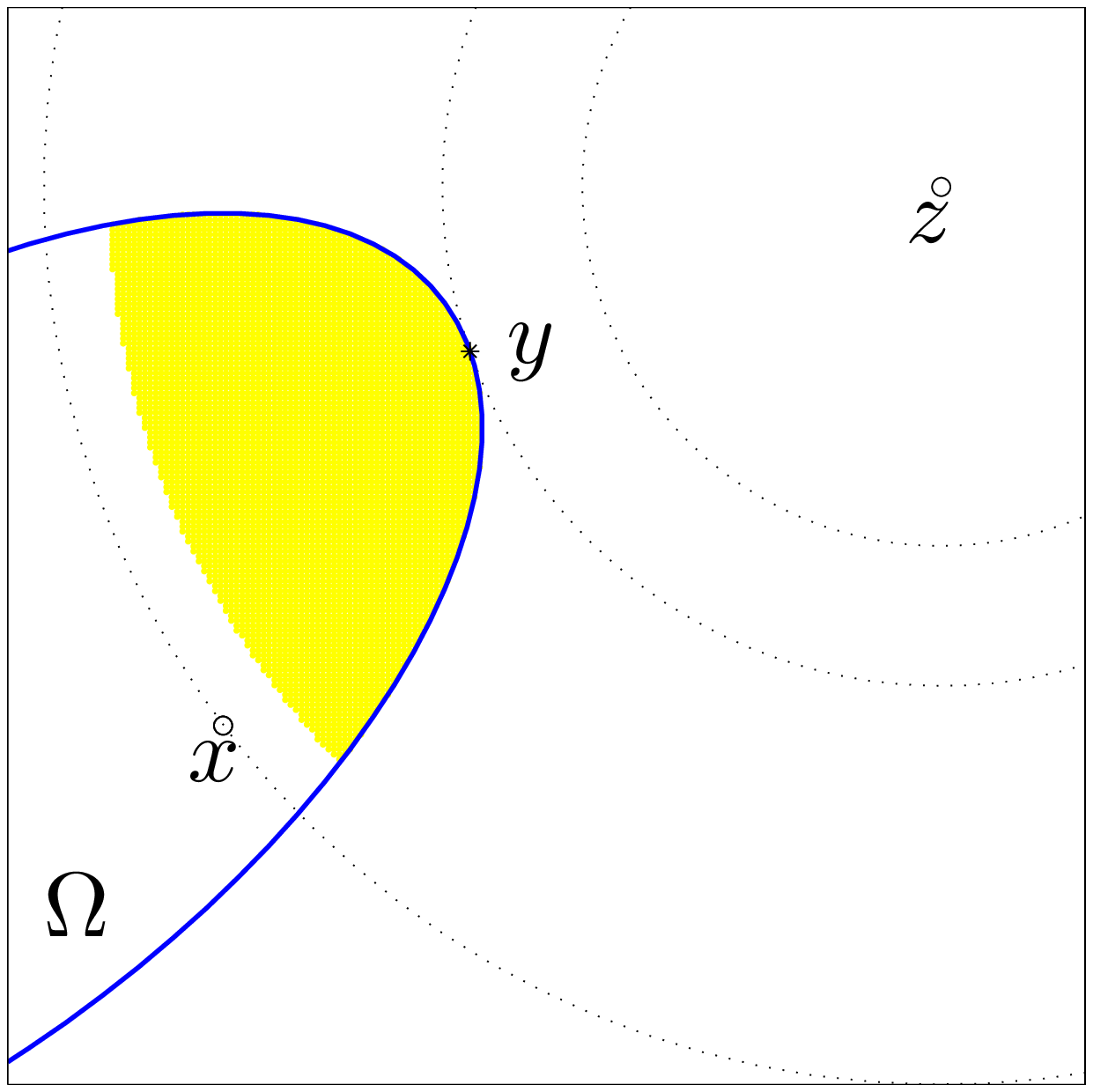}}&{\includegraphics[scale=\FigScala]{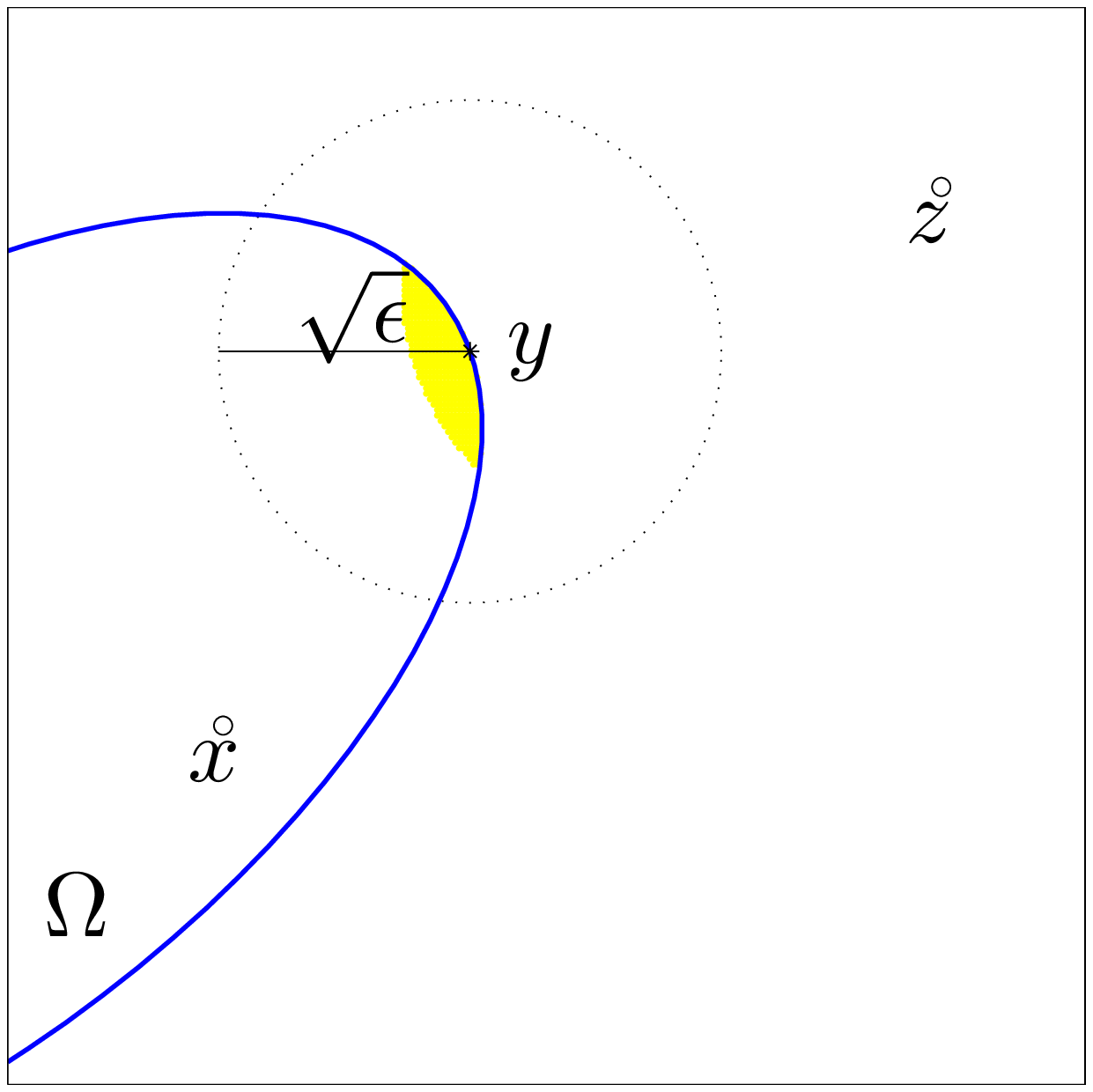}}\\
\end{tabular}
\end{center}
\caption{Example with $f_1(\x) = \iota_\Omega(\x)$, $\ds$ as in \eqref{scaled-euclidean} with $\alpha = 1$, $D=I$. Left panel: in yellow, the set $P_\eta(\x;\hs)$ defined in \eqref{etaset}. Right panel: in yellow, the set of points $\ty$ satisfying \eqref{eps-approx}.}\label{fig:inexact-pictures}
\end{figure}

In general, one of the main differences between definitions \eqref{eta-approxk} and \eqref{eps-approxk} consists in the fact that in the latter case the distance between the approximated and the exact minimum of $\hskxy{\puntino}\xk$, i.e. $\|\tyk-\yk\|$, can be controlled by the independent parameter $\epsilon_k$, while in the other case this distance is algorithm and iteration dependent. This fact can be exploited to obtain a stronger convergence result, as shown in the next section.

\subsection{Convergence analysis in the convex case with $\epsilon$-approximations}\label{sec:convergence}

\subsubsection{Convergence}

In this section, we assume that $f_0$ is convex and, in this case, we prove a stronger convergence result for a specific line--search algorithm where the descent direction is defined by means of an $\epsilon$-approximation, provided that the sequence of parameters $\{\epsilon_k\}_\kinN$ is summable and that the sequence of the matrices $D_k$ satisfies suitable assumptions. The following theorem is a generalization of Theorem 3.1 in \cite{Bonettini-Prato-2015a}. Further results on forward-backward variable metric algorithms which apply to problems of the form \eqref{minf} when $f_0$ has Lipschitz continuous gradient can be found in the recent papers \cite{Combettes-Vu-2014,Combettes-Vu-2013}. We stress that in all our analysis we do not need any Lipschitz continuity of the gradient of $f_0$ and, moreover, the sequence of errors $\|\tyk-\yk\|$ needs to be square summable, while the convergence result stated in \cite[Theorem 4.1]{Combettes-Vu-2014} is given under the stronger assumption that $\|\tyk-\yk\|$ is summable.
\begin{Thm}\label{thm:sgp_convex_converge}
Let $0<\alpha_{\min}\leq \alpha_{\max}$, $\gamma \in [0,1]$, $\{\ak\}_\kinN\subset[\alpha_{\min},\alpha_{\max}]$. Assume that $f_0$ in \eqref{minf} is convex and the solution set $X^*$ of problem \eqref{minf} is not empty. Let $\{\xk\}_{k \in \N}$ be the sequence generated as
\begin{equation*}
\xkk = \xk + \lamk\dk,\ \ \ \dk = \tyk-\xk
\end{equation*}
where $\lamk$ is obtained by means of the backtracking procedure in Algorithm \ref{Algo2}, with $\tyk$ satisfying $\thskxy\tyk\xk<0$. Moreover assume that:
\begin{enumerate}[label=\textup{(H{\arabic*})}]
\item \label{hypot1} $\tyk$ satisfies \eqref{eps-approxk}, where the sequence $\{\epsilon_k\}_{\kinN}$ is summable, i.e. $\sum_{k=0}^{\infty} \epsilon_k < \infty$;
\item \label{hypot2} $\{D_k\}_\kinN\subset {\mathcal M}_\mu$, where $\mu \geq 1$ and
    \begin{equation}\nonumber
    D_{k+1} \preceq (1+\zeta_k)D_k, \ \ \ \{\zeta_k\}_\kinN \subset \R_{\geq0},\quad\mathrm{and}\quad \sum_{k=0}^\infty \zeta_k < \infty.
    \end{equation}
\end{enumerate}
Then the sequence $\{\xk\}_{k \in \N}$ converges to a solution of \eqref{minf}.
\end{Thm} \noindent
\proof First of all we recall the basic norm equality
\begin{equation}\label{3pointformula}
\|a-b\|_D^2 +\|b-c\|_D^2-\|a-c\|^2_D = 2(a-b)^TD(c-b)
\end{equation}
which holds for any $a,b,c\in \R^n$. Let $\hatx \in X^*$. By definition of $\tyk$ we have
\begin{equation*}
f_1(\w)\geq f_1(\tyk)+\frac 1{\ak}(\zk-\tyk)^TD_k(\w-\tyk) - \epsilon_k\ \ \ \forall \w\in\R^n
\end{equation*}
which, recalling that $\zk = \xk-\ak D_k^{-1}\nabla f_0(\xk)$, writes also as
\begin{equation*}
(\tyk-\xk)^TD_k(\w-\tyk)\geq\ak \left(f_1(\tyk)-f_1(\w)+\nabla f_0(\xk)^T(\tyk-\w)\right)-\ak\epsilon_k \ \ \forall\w\in\R^n.
\end{equation*}
For $\w = \hatx$, the previous inequality gives
\begin{align}
(\tyk-\xk)^TD_k(\hatx-\xk) &\geq \ak\left(f_1(\tyk)-f_1(\hatx) +\nabla f_0(\xk)^T(\xk-\hatx)\right)-\ak\epsilon_k+\nonumber\\
   &  \qquad\qquad+\left(\tyk-\xk+\ak D_k^{-1}\nabla f_0(\xk)\right)^TD_k(\tyk-\xk)\nonumber\\
&\geq \ak\left(f_1(\tyk)-f_1(\xk)+f(\xk)-f(\hatx)\right)+ \|\tyk-\xk\|^2_{D_k}\label{add2}\\
   & \qquad+\ak\nabla f_0(\xk)^T(\tyk-\xk)-\ak\epsilon_k\nonumber\\
&\geq \|\tyk-\xk\|^2_{D_k}-\ak\epsilon_k\nonumber\\
   & \qquad\qquad+\ak\left(f_1(\tyk)-f_1(\xk)+\nabla f_0(\xk)^T(\tyk-\xk)\right)\nonumber\\
   &= \frac{1}{(\lamk)^2}\|\xkk-\xk\|_{D_k}^2-\ak\epsilon_k\label{tmp1}\\
   & \qquad\qquad+ \ak\left(f_1(\tyk)-f_1(\xk)+\nabla f_0(\xk)^T(\tyk-\xk)\right),\nonumber
\end{align}
where the second inequality is obtained adding and subtracting $f_1(\xk)$ and by the convexity of $f_0$, the third one from the fact that $\hatx$ is a minimum point and the last one by definition of $\xkk$. By equality \eqref{3pointformula} with $a =\xkk$, $b = \xk$, $c = \hatx$, $D=D_k$ we obtain
\begin{eqnarray}
\|\xkk-\hatx\|^2_\Dk&=& \|\xk-\hatx\|^2_\Dk+\|\xkk-\xk\|^2_\Dk-2(\xk-\xkk)^T\Dk(\xk-\hatx)\nonumber\\
&=&\|\xk-\hatx\|^2_\Dk+\|\xkk-\xk\|^2_\Dk-2\lamk(\tyk-\xk)^T\Dk(\hatx-\xk)\nonumber\\
&\stackrel{\eqref{tmp1}}{\leq}&\|\xk-\hatx\|^2_\Dk+\left(1-\frac 2 {\lamk}\right)\|\xkk-\xk\|^2_\Dk+\nonumber\\
& &-2\ak\lamk\left(\nabla f_0(\xk)^T(\tyk-\xk)+f_1(\tyk)-f_1(\xk)\right)+2\ak\lamk\epsilon_k\nonumber\\
&=&\|\xk-\hatx\|^2_\Dk+\left(1-\frac 2 {\lamk}+\frac \gamma {\lamk}\right)\|\xkk-\xk\|^2_\Dk+\nonumber\\
& &-2\ak\lamk\thskxy\tyk\xk+2\ak\lamk\epsilon_k\nonumber\\
&\leq&\|\xk-\hatx\|^2_\Dk-2\ak\lamk\thskxy\tyk\xk+2\ak\lamk\epsilon_k,\label{add3}
\end{eqnarray}
where the third equality is obtained by adding and subtracting the term $\gamma\lamk\|\tyk-\xk\|_\Dk^2=\gamma/\lamk \|\xkk-\xk\|_\Dk^2$ and the last inequality follows from the fact that $\gamma \in [0,1]$.
From assumption \ref{hypot2} we obtain
\begin{eqnarray}
\|\xkk-\hatx\|^2_\Dkk &\leq & (1+\zeta_k)\|\xkk-\hatx\|^2_\Dk\nonumber\\ &\leq&(1+\zeta_k)\|\xk-\hatx\|^2_\Dk-2\ak(1+\zeta_k)\lamk\thskxy\tyk\xk\nonumber\\
&& \qquad\qquad+2\ak\lamk(1+\zeta_k)\epsilon_k\nonumber\\
&\leq&(1+\zeta_k)\|\xk-\hatx\|^2_\Dk-2\alpha_{\max}\zeta\lamk\thskxy\tyk\xk+2\alpha_{\max}\zeta\epsilon_k\label{add6}
\end{eqnarray}
where we set $\zeta = 1+ \max_k \zeta_k$. Then, from \cite[Lemma 2.2.2]{Polyak-1987} we can conclude that the sequence $\{\|\xk-\hatx\|^2_\Dk\}_\kinN$ converges. In particular, since $D_k\in {\mathcal M}_\mu$, $\{\xk\}_\kinN$ is bounded and, thus, it has at least one limit point. Let us denote such limit point by $\x^{\infty}$. By Corollary \ref{cor:eps-approx}, $\x^{\infty}$ is stationary; in particular, since $f$ is convex, it is a minimum point, i.e. $\x^{\infty}\in X^*$ and, thus, $\{\|\xk-\x^\infty\|^2_\Dk\}_\kinN$ converges. Let $\{\x^{(k_i)}\}_{i \in \N}$ be a subsequence of $\{\xk\}_{k \in \N}$ which converges to $\x^\infty$. By the norm inequality \eqref{ine_norm} we can write
\begin{equation}\nonumber
\|\x^{(k_i)}-\x^\infty\|_{D_{k_i}}^2\leq \mu \|\x^{(k_i)}-\x^\infty\| \stackrel{i\to\infty}{\longrightarrow} 0
\end{equation}
Since $\{\|\xk-\x^\infty\|^2_\Dk\}_\kinN$ converges, this implies that its limit is zero. Invoking again \eqref{ine_norm} we can write
\begin{equation}\nonumber
\frac 1 \mu \|\x^{(k)}-\x^\infty\|^2\leq\|\x^{(k)}-\x^\infty\|_{D_{k}}^2 \stackrel{k\to\infty}{\longrightarrow} 0
\end{equation}
which allows to conclude that $\{\xk\}_\kinN$ converges to $\x^{\infty}$.
\endproof
In the following we present a variation of Theorem \ref{thm:sgp_convex_converge} where the tolerance parameters $\ek$ are adaptively chosen, instead of being a prefixed summable sequence.
\begin{Thm}\label{thm:add1}
Let $0<\alpha_{\min}\leq \alpha_{\max}$, $\gamma \in [0,1]$, $\{\ak\}_\kinN\subset[\alpha_{\min},\alpha_{\max}]$. Assume that $f_0$ in \eqref{minf} is convex and the solution set $X^*$ of problem \eqref{minf} is not empty. Let $\{\xk\}_{k \in \N}$ be the sequence generated as
\begin{equation*}
\xkk = \xk + \lamk\dk,\ \ \ \dk = \tyk-\xk
\end{equation*}
where $\lamk$ is obtained by means of the backtracking procedure in Algorithm \ref{Algo2}, with $\tyk$ satisfying $\thskxy\tyk\xk<0$. Moreover assume that
\begin{enumerate}[label=\textup{(H{\arabic*}')}]
\item\label{hypot1prime} $\tyk$ satisfies \eqref{eps-approxk}, where the sequence $\{\epsilon_k\}_{\kinN}$ satisfies
    \begin{equation}\label{eps-tau}
    \ek \leq -\tau \thskxy\tyk\xk
    \end{equation}
    for some $\tau > 0$,
\end{enumerate}
and that hypothesis \ref{hypot2} of Theorem \ref{thm:sgp_convex_converge} holds. Then, the sequence $\{\xk\}_{k \in \N}$ converges to a solution of \eqref{minf}.
\end{Thm}
\proof By substituting \eqref{eps-tau} in \eqref{add6} we obtain
\begin{equation}\nonumber
\|\xkk-\hatx\|^2_\Dkk\leq(1+\zeta_k)\|\xk-\hatx\|^2_\Dk-2\alpha_{\max}\zeta(1+\tau)\lamk\thskxy\tyk\xk.
\end{equation}
The rest of the proof follows exactly from the same arguments employed in Theorem \ref{thm:sgp_convex_converge}.\endproof
We will show in Section~\ref{sec:eps-implementation} how the conditions \eqref{eps-approxk} and \eqref{eps-tau} can be satisfied in practice.

Assumption \ref{hypot2} is analogous to the one proposed in \cite{Combettes-Vu-2014,Combettes-Vu-2013}. A special case of it consists in the following
\begin{enumerate}[label=\textup{(H{\arabic*}')}]\setcounter{enumi}{1}
\item\label{hypot2prime} $\{\Dk\}_\kinN\subset{\mathcal M}_{\mu_k}$, where
$\displaystyle    \mu_k^2 = 1+\xi_k,\ \ \xi_k\geq 0, \ \ \ \sum_{k=0}^{\infty}\xi_k<\infty.$
\end{enumerate}
Thanks to the inequality \eqref{ine_norm}, for any $\x\in\R^n$ we have
\begin{equation*}
\x^T(\Dkk-\mu_k\mu_{k+1}D_k)\x = \x^T\Dkk\x-\mu_k\mu_{k+1}\x^TD_k\x \leq \mu_{k+1}\|\x\|^2 - \mu_k\mu_{k+1}\frac{\|\x\|^2}{\mu_k} = 0,
\end{equation*}
which implies $\Dkk \preceq \mu_k\mu_{k+1}D_k$. Moreover, $\mu_k\mu_{k+1}$ can be written as $\mu_k\mu_{k+1} = 1+ \zeta_k$, where $\zeta_k = \sqrt{(1+\xi_k)(1+\xi_{k+1})}-1$. Since $\lim_{x\to 0}\sqrt{1+x}/x = 1/2$, it follows that $\sum_{k=0}^\infty \xi_k$ and $\sum_{k=0}^\infty \zeta_k$ have the same behaviour. Then, we can conclude that (H2') implies (H2).

We also observe that, employing the same arguments above, we can also prove that $\mu_{k+1}\mu_k\Dkk \succeq \Dk$, and, as a consequence, \ref{hypot2prime} also implies that $(1+\zeta_k) \Dkk \succeq \Dk$ with $\sum_{k=0}^\infty \zeta_k<\infty$.\\
In practice, \ref{hypot2prime} says that the scaling matrices have to converge to the identity matrix at a certain rate, while \ref{hypot2} implies the convergence to some symmetric positive definite matrix (see Lemma 2.3 in \cite{Combettes-Vu-2013}).

\subsubsection{Convergence rate analysis}

In this section we analyze the convergence rate of the objective function values $f(\xk)$ to the optimal one, $f^*$, proving that $f(\xkk)-f^* = {\mathcal O}(\frac 1 k)$. This complexity result is obtained in the same settings of Theorem \ref{thm:add1}, but further assuming that the gradient of $f_0$ is Lipschitz continuous on the domain of $f_1$. This Lipschitz assumption guarantees that the sequence $\{\lamk\}_\kinN$ is bounded away from zero.
Before giving the main results, we need to prove the following lemma, which actually does not require the Lipschitz assumption.
\begin{Lemma}\label{lemma:1_new}
Let $\xk,\tyk\in \Omega$. If $\tyk$ satisfies \eqref{eps-approxk}, with $0<\ak\leq\alpha_{\max}$ and $D_k \in {\mathcal M}_\mu$, then,
\begin{equation}\label{new_work2}
\frac{1}{2\alpha_{\max}\mu}\|\tyk-\xk\|^2 \leq  -\thskxy\tyk\xk +\epsilon_k.
\end{equation}
\end{Lemma}
\proof For any $\w\in \partial_{\epsilon_k} f_1(\tyk)$ we have
\begin{eqnarray*}
\hskxy\tyk\xk &=& \nabla f_0(\xk)^T(\tyk-\xk) + \frac 1{2\ak}\|\tyk-\xk\|^2_{D_k} + f_1(\tyk)-f_1(\xk)\\
&\leq & \nabla f_0(\xk)^T(\tyk-\xk) + \frac 1{2\ak}\|\tyk-\xk\|^2_{D_k} + \w^T(\tyk-\xk) + \epsilon_k.
\end{eqnarray*}
In particular, the previous inequality holds true for $\w = \frac{1}{\ak} D_k(\zk-\tyk)$ (see \eqref{eps-approxk}). This results in
\begin{eqnarray}
\thskxy\tyk\xk &\leq & \hskxy\tyk\xk\nonumber\\
&\leq&\nabla f_0(\xk)^T(\tyk-\xk) + \frac 1{2\ak}\|\tyk-\xk\|^2_{D_k} + \nonumber\\
& & + \frac{1}{\ak} (\xk-\ak D_k^{-1} \nabla f_0(\xk)-\tyk)^TD_k(\tyk-\xk) + \epsilon_k\nonumber\\
&=& -\frac{1}{2\ak}\|\tyk-\xk\|^2_{D_k} +\epsilon_k\leq -\frac{1}{2\alpha_{\max}\mu}\|\tyk-\xk\|^2 +\epsilon_k,\nonumber
\end{eqnarray}
where the last inequality follows from \eqref{ine_norm}.
\hfill$\square$
\begin{Pro}\label{pro:Lipschitz}
Let $\{\xk\}_{k\in\N}$ be a sequence of points in $\Omega$ and $\{\dk\}_\kinN$ a sequence of descent directions such that $\dk = \tyk-\xk$ and \eqref{new_work2} holds. Let $\{\lamk\}_{k\in\N}$ be the steplength sequence computed by Algorithm \ref{Algo2} and assume that $\nabla f_0$ is Lipschitz continuous on $\Omega$ and that \eqref{eps-tau} holds.
Then, there exists  $\lambda_{\min}\in \R_{>0}$ such that
\begin{equation}\label{lambdastar}
\lamk\geq \lambda_{\min} \ \ \ \forall k\in\N. 
\end{equation}
\end{Pro}
\proof In view of \eqref{eps-tau}--\eqref{new_work2}, setting $a = \alpha_{\max}\mu$, one obtains
\begin{equation}\label{inpro}
\|\dk\|^2\leq -2a(1+\tau)\thskxy\tyk\xk.
\end{equation}
If $\nabla f_0$ is Lipschitz continuous on $\Omega$ with Lipschitz constant $L$, then from the descent lemma \cite[p.667]{Bertsekas-1999} we have
\begin{equation}\label{descent-lemma}
f_0(\xk+\lambda \dk) \leq f_0(\xk) +\lambda \nabla f_0(\xk)^T\dk +\frac L 2 \lambda^2\|\dk\|^2,
\end{equation}
where $\lambda\in[0,1]$.
By combining inequalities \eqref{inpro} and \eqref{descent-lemma} we further obtain
\begin{equation}\nonumber
f_0(\xk+\lambda \dk) \leq f_0(\xk) +\lambda\nabla f_0(\xk)^T\dk-a(1+\tau)L  \lambda^2\thskxy\tyk\xk.
\end{equation}
Summing $f_1(\xk+\lambda\dk)$ on both sides of the previous relation and applying the Jensen inequality $f_1(\xk+\lambda\dk)\leq (1-\lambda)f_1(\xk)+\lambda f_1(\tyk)$ to the r.h.s. yields
\begin{eqnarray*}
f(\xk+\lambda \dk) &\leq& f(\xk) -\lambda f_1(\xk)+\lambda f_1(\tyk) + \lambda\nabla f_0(\xk)^T\dk\\
&&\qquad-aL  \lambda^2(1+\tau)\thskxy\tyk\xk\\
&\leq& f(\xk) - \lambda f_1(\xk)+\lambda f_1(\tyk)+\lambda\nabla f_0(\xk)^T\dk\\
& &\qquad -aL  \lambda^2(1+\tau)\thskxy\tyk\xk + \frac{\lambda \gamma}2\|\dk\|^2_{D_k}\\
&=& f(\xk) +\lambda \thskxy\tyk\xk -aL  \lambda^2(1+\tau)\thskxy\tyk\xk\\
&=& f(\xk) +\lambda\left(1-aL(1+\tau)\lambda\right) \thskxy\tyk\xk.
\end{eqnarray*}
The previous inequality ensures that the Armijo condition
\begin{equation}\label{armijo2}
f(\xk+\lambda\dk)\leq f(\xk)+\lambda\beta\thsxy\tyk\xk\end{equation}
is satisfied, for all $k\in\N$, when $1-aL(1+\tau)\lambda\geq \beta$, that is for all $\lambda$ such that $\lambda \leq (1-\beta)/(aL(1+\tau))$. If $\lamk$ is the steplength computed by Algorithm \ref{Algo2} and the backtracking loop is performed at least once, then $\lambda =\lamk/\delta$ does not satisfy inequality \eqref{armijo2}, which means $\lamk > (1-\beta)\delta/(aL(1+\tau))$. Thus, the steplength sequence $\{\lamk\}_{k\in\N}$ satisfies inequality \eqref{lambdastar} with $\lambda_{\min} = (1-\beta)\delta/(aL(1+\tau))$.
\endproof
Based on these premises, we are now ready to prove the convergence rate result.
\begin{Thm}\label{thm:SGP_convergence_rate}
Assume that the hypotheses of Theorem \ref{thm:add1} hold and, in addition, that the gradient of $f_0$ is Lipschitz continuous on $\Omega$. Let $f^*$ be the optimal function value for problem \eqref{minf}. Then, we have
\begin{equation*}
f(\xkk)-f^* = {\mathcal O} \left(\frac 1 k\right).
\end{equation*}
\end{Thm}
\proof
If we do not neglect the term $f(\xk)-f(\hatx) = f(\xk)-f^*$ in \eqref{add2} and in all the subsequent inequalities, instead of \eqref{add3} we obtain
\begin{eqnarray*}
\|\xkk-\hatx\|^2_\Dk &\leq&\|\xk-\hatx\|^2_\Dk+2\ak\lamk\left(-\thskxy\tyk\xk+\epsilon_k\right)-2\lamk\ak(f(\xk)-f^*),
\end{eqnarray*}
and hence:
\begin{eqnarray*}
\|\xkk-\hatx\|^2_\Dkk &\leq& (1+\zeta_k)\|\xkk-\hatx\|^2_\Dk\\ &\leq&(1+\zeta_k)\|\xk-\hatx\|^2_\Dk+2\ak\lamk(1+\zeta_k)(-\thskxy\tyk\xk+\epsilon_k) +\\
& &\qquad -2\lamk(1+\zeta_k)\ak(f(\xk)-f^*)\\
&\stackrel{\eqref{eps-tau}}{\leq}&(1+\zeta_k)\|\xk-\hatx\|^2_\Dk-2\alpha_{\max}(1+\tau)\zeta\lamk\thskxy\tyk\xk +a(f^*-f(\xk)),
\end{eqnarray*}
where we set $\zeta = 1+ \max_k \zeta_k$, $a = 2\lambda_{\min}\alpha_{\min}$, where $\lambda_{\min}$ is defined in Proposition \ref{pro:Lipschitz}.
Summing the previous inequality from 0 to $k$ gives
\begin{eqnarray*}
\|\xkk-\hatx\|^2_\Dkk &\leq&\|\x^{(0)}-\hatx\|^2_{D_0} + \sum_{j=0}^k \zeta_j \|\xj-\hatx\|_{D_j} ^2 -2\alpha_{\max}(1+\tau)\zeta\sum_{j=0}^k\lamj\thsjxy\tyj\xj+\\ & & +a\left((k+1)f^*-\sum_{j=0}^kf(\xj)\right)\\
&\leq&\|\x^{(0)}-\hatx\|^2_{D_0} + M\bar\zeta -\frac{2\alpha_{\max}(1+\tau)\zeta}\beta(f(\x^{(0)})-f^*) +a\left((k+1)f^*-\sum_{j=0}^kf(\xj)\right),
\end{eqnarray*}
where the second inequality follows by setting $\bar\zeta = \sum_{j=0}^\infty\zeta_j$, from the fact that $\{\|\xk-\hatx\|_{D_k}^2\}_\kinN$ is a convergent sequence (see Theorem \ref{thm:add1}), thus there exists $M$ such that $\|\xj-\hatx\|_{D_j}^2\leq M$, and from \eqref{add4}.
Adding the positive quantity $a(f(\x^{(0)}) - f^*)$ to the right hand side of the last inequality we obtain
\begin{equation*}
\|\xkk-\hatx\|^2_\Dkk
\leq\|\x^{(0)}-\hatx\|^2_{D_0} + M\bar\zeta -\frac{2\alpha_{\max}(1+\tau)\zeta}\beta(f(\x^{(0)})-f^*) +a\left(kf^*-\sum_{j=1}^kf(\xj)\right).
\end{equation*}
Moreover, exploiting the inequality
\begin{equation}\nonumber
0\leq \sum_{j=0}^k j (f(\xj)-f(\xjj)) = \sum_{j=1}^k f(\xj) - kf(\xkk)
\end{equation}
gives
\begin{equation*}
\|\xkk-\hatx\|^2_\Dkk
\leq\|\x^{(0)}-\hatx\|^2_{D_0} + M\bar\zeta -\frac{2\alpha_{\max}(1+\tau)\zeta}\beta(f(\x^{(0)})-f^*) +ak(f^*-f(\xkk)).
\end{equation*}
Rearranging terms, this finally yields
\begin{equation*}
f(\xkk) - f(\hatx)\leq \frac 1{ak}\left( \|\x^{(0)}-\hatx\|^2_{D_0} + M\bar\zeta - 2\frac{\alpha_{\max}(1+\tau)\zeta}{\beta}(f(\x^{(0)})-f(\hatx))\right),
\end{equation*}
establishing the result.\hfill $\square$

\section{Practical computation of $\eta$- and $\epsilon$- approximations}
\subsection{Computing $\eta$-approximations}\label{sec:eta-implementation}
In this section we discuss how to compute a point $\tyk$ such that \eqref{eta-approxk} holds, i.e. satisfying
\begin{equation}\label{eta-approx}
\hskxy\tyk\xk\leq\eta \hskxy\yk\xk, \ \ \mbox{ with } \yk = p(\xk;\hs),
\end{equation}
for a given $\eta \in (0,1]$, without knowing $\yk$. A special case of this problem, corresponding to the case $f_1=\iota_\Omega$, where $\Omega$ is the intersection of closed, convex sets and the metric is given by \eqref{scaled-euclidean}, is considered in \cite{Birgin-etal-2003}. This is possible when, for each $k$, one can compute a sequence $\{a_l\}_{l\in\N}\subset \R$ such that
\begin{equation}\label{aprop}
a_l\leq \hskxy\yk\xk, \ \ \forall l\in \N, \ \ \mbox{ and } \lim_{l\rightarrow \infty} a_l = \hsxy\yk\xk,
\end{equation}
and a sequence of points $\{\tykl\}_{l\in\N}$ such that
\begin{equation}\label{aprop2}
\lim_{l\rightarrow \infty} \hskxy\tykl\xk = \hskxy\yk\xk.
\end{equation}
In practice, $l$ should be considered as the index of an inner loop for computing $\tyk$. Indeed, when \eqref{aprop} holds, we also have
\begin{equation}\label{aprop1}
\eta a_l \leq \eta \hskxy\yk\xk \ \ \forall l\in \N.
\end{equation}
Moreover, for all sufficiently large $l$ we have $a_l> \hskxy\yk\xk/\eta$ which, together with \eqref{aprop1} gives
\begin{equation*}
\hskxy\yk\xk < \eta a_l \leq \eta \hskxy\yk\xk.
\end{equation*}
Then, if one considers any method generating a sequence $\tykl$ such that \eqref{aprop2} holds, the stopping criterion
\begin{equation}\label{stoprule}
\hskxy\tykl\xk \leq \eta a_l
\end{equation}
for the inner iterations is well defined. If $l$ is the smallest integer such that \eqref{stoprule} is satisfied, then the point $\tyk=\tykl$ satisfies \eqref{eta-approx}. In the following sections we show how to compute a sequence $a_l$ satisfying \eqref{aprop} in an interesting case.

\subsection{Composition with a linear operator}\label{sec:affine}

In this section we assume that $f_1(x)$ is given by
\begin{equation}\label{affine}
f_1(\x) = g(A\x),
\end{equation}
where $A\in \R^{m\times n}$ and $g: \R^m\rightarrow \bR$ is a convex function. Moreover, we choose $\ds$ as in \eqref{scaled-euclidean}. Let us consider the minimum problem \eqref{proj} which can be written in equivalent primal--dual and dual form as
\begin{equation*}
\min_{y\in\R^n}\hskxy\y\xk = \min_{y\in\R^n}\max_{\v\in \R^m} \Fsk(\y,\v,\xk) = \max_{\v\in \R^m} \Psisk(\v,\xk).
\end{equation*}
The primal--dual problem can be obtained from the primal one by applying Definition \ref{def:dual} of the convex conjugate, which gives $g(A\x) = \max_{\v\in\R^m} \v^TA\x-g^*(v)$, obtaining
\begin{equation}\label{primal-dual-function}
\Fsk(\y,\v,\xk)=\frac 1{2\ak}\|y-\zk\|^2_{D_k} + \y^TA^T\v - g^*(\v)-f_1(\xk)-\frac{\ak}{2} \|\nabla f_0(\xk)\|^2_{D_k^{-1}}
\end{equation}
with $\zk = \xk-\ak D_k^{-1}\nabla f_0(\xk)$. The dual problem is obtained by computing the minimum of the primal--dual function with respect to $\y$, which is given by  $\y=\zk-\ak D_k^{-1} A^T v$, and substituting it in \eqref{primal-dual-function}, obtaining the explicit expression of the dual function
\begin{equation*}
\Psiskxy\v = -\frac{1}{2\ak}\|\ak D_k^{-1}A^T\v-\zk\|^2_{D_k} - g^*(\v)-f_1(\xk)-\frac{\ak}{2} \|\nabla f_0(\xk)\|^2_{D_k^{-1}} + \frac{1}{2\ak}\|\zk\|^2_{D_k}.
\end{equation*}
By definition of the primal--dual and dual functions, the following inequalities hold
\begin{equation}\nonumber
\hskxy\y\xk\geq \Fsk(\y,\v,\xk)\geq \Psisk(\v,\xk) \ \ \forall \y\in \R^n,\v\in\R^m.
\end{equation}
In particular, the previous inequality holds for $\y = \yk$.
Then, an approximation $\tyk$ of $\yk$ can be computed by applying any method to the dual problem
\begin{equation}\label{dual}
\max_{v\in\R^m}  \Psisk(\v,\xk),
 \end{equation}
generating a sequence $\{\vl\}_{l\in\N}$ such that $\Psiskxy\vl$ converges to the maximum of the dual function $\Psiskxy{\puntino}$. As a consequence of this, setting $\tykl = \zk -\ak D_k^{-1} A^T \vl$, a point satisfying \eqref{eta-approx} can be found by stopping the dual iterations when
\begin{equation}\label{stoprule1}
\hskxy\tykl\xk\leq\eta\Psisk(\vl,\xk)
\end{equation}
is satisfied, i.e. \eqref{stoprule} with $a_l =  \Psiskxy\vl$.

For example, one can apply a forward--backward method \cite{Combettes-Pesquet-2011}, called also ISTA or its accelerated version (FISTA, \cite{Beck-Teboulle-2009b}) to the dual problem. As an alternative, also the saddle point problem
\begin{equation}\nonumber
\min_{\y\in \R^n }\max_{v\in\R^m} \Fsk(\y,\v,\xk)
 \end{equation}
can be faced, for example with a primal--dual method such as \cite{Chambolle-Pock-2011,Loris-Verhoeven-2011}, using \eqref{stoprule1} as stopping condition.
More in general, a point $\tyk\in P_\eta(\xk;\hsk)$ can be obtained by
computing two sequences, $\{\vl\}_{l\in\N}$, $\{\tykl\}_{l\in\N}$, such that
\begin{equation}\nonumber
\lim_{l\rightarrow\infty}\Psiskxy\vl= \max_{v\in\R^m} \Psiskxy\v= \min_{\y\in\R^n} \hskxy\y\xk=\lim_{l\rightarrow\infty}\hskxy\tykl\xk,
\end{equation}
stopping the iterates when \eqref{stoprule1} is met.
\paragraph{Remarks}
We observe that \eqref{affine} includes also the case where $f_1(\x)$ is defined as $f_1(x) = \sum_{i=1}^r g_i(A_i\x)$, where $A_i\in \R^{m_i\times n}$, $g_i:\R^{m_i}\rightarrow \R$. Indeed, formulation \eqref{affine} is recovered by setting $A = [A_1^T\ A_2^T\ ... \ A_r^T]^T\in \R^{m\times n}$ with $m = \sum_{i = 1}^r m_i$. In this case the dual variable $\v$ can be partitioned as $\v = [\v_1^T\ \v_2^T ... \ \v_r^T]^T$, where $\v_i\in \R^{m_i}$ and $g^*(v) = \sum_{i=1}^r g_i^*(\v_i)$ (see \cite[Theorem 2.3.1 (iv)]{Zalinescu-2002}).  

\subsection{Preserving feasibility}

Clearly, any point $\tykl$ satisfying \eqref{stoprule1}, where $\vl$ is generated by any converging algorithm applied to the dual or the primal--dual problem, belongs to the domain of $\hs(\cdot,\xk)$, i.e. to the set $\Omega$. Indeed, for any $l$, $\vl$ belongs to the domain of the dual function $\Psisk(\cdot,\xk)$ and, as a consequence, \eqref{stoprule} implies that $\hskxy\tykl\xk$ is finite. However, the stopping criterion \eqref{stoprule} may require a very large number of inner iterations $l$ to be satisfied, and, in addition, the primal sequence points $\tykl$ may be feasible only in the limit. For these reasons, we propose to consider also the sequence $\bykl = P_\Omega(\tykl)$, where $P_\Omega$ denotes the Euclidean projection onto the set $\Omega$. If, at some inner iteration $l$, the inequality
\begin{equation}\label{stoprule2}
\bthskxy\bykl\xk \leq \eta \Psiskxy\vl
\end{equation}
is satisfied, this clearly means that $\bykl\in P_\eta(\xk;\hs)$ (i.e., \eqref{eta-approxk} is satisfied) and we can set $\tyk = \bykl$. We observe that, when $\tykl$ converges to $\yk$ as $l$ diverges, the stopping criterion \eqref{stoprule2} is well defined, since $\bykl$ also converges to $\yk$.

\subsection{Computing $\epsilon$-approximations}\label{sec:eps-implementation}

In this section we show how to compute a point satisfying inclusion \eqref{eps-approxk}, for any given $\epsilon_k \in\R_{\geq 0}$, when the convex function $f_1$ in \eqref{minf} has the form \eqref{affine}. Our arguments are obtained by extending those in \cite{Villa-etal-2013}, which are recovered setting $D_k = I$. As done in Section \ref{sec:affine}, we will make use of the duality theory. In particular, we define the primal--dual gap function as
\begin{equation}\label{gap}
\Gskxy\y\v = \hskxy\y\xk - \Psiskxy\v.
\end{equation}
We also have the following result.
\begin{Pro}
Let $\epsilon_k\in\R_{\geq 0}$. If
\begin{equation}\label{gap-cond}
\Gskxy\tyk\v\leq \epsilon_k,
\end{equation}
with $\tyk = \zk-\ak D_k^{-1}A^Tv$, for some $v\in \R^m$, then \eqref{eps-approxk} is satisfied.
\end{Pro}
\proof From the definition of the primal--dual gap, a simple computation shows that
\begin{eqnarray*}
\Gskxy\tyk\v &=& \frac{1}{\ak}\|\ak D_k^{-1}A^T\v\|^2_{D_k} - \v^TA\zk + f_1(\zk-\ak D_k^{-1}A^T\v) + g^*(\v)\\
&=&\sup_{w\in\R^m} \frac{1}{\ak}\|\ak D_k^{-1}A^T\v\|^2_{D_k} - \v^TA\zk + \w^T(\zk-\ak D_k^{-1}A^T\v)-f_1^*(\w)+g^*(\v)\\
&=&\sup_{w\in\R^m} (w-A^Tv)^T(\zk-\ak D_k^{-1}A^T\v) -f_1^*(\w)+g^*(\v)\\
&\geq& \sup_{w\in\R^m} (w-A^Tv)^T(\zk-\ak D_k^{-1}A^T\v) -f_1^*(\w)+f_1^*(A^T\v),
\end{eqnarray*}
where the last inequality follows from Proposition \ref{pro:duality}. Thus, if \eqref{gap-cond} holds, the previous inequality yields
\begin{equation*}
(w-A^Tv)^T(\zk-\ak D_k^{-1}A^T\v) -f_1^*(\w)+f_1^*(A^T\v) \leq \Gskxy\tyk\v\leq\epsilon_k \ \ \ \forall w\in \R^m.
\end{equation*}
Rearranging terms, the previous inequality writes also as
\begin{equation*}
f_1^*(\w)\geq f_1^*(A^T\v) +(w-A^Tv)^T(\zk-\ak D_k^{-1}A^T\v) -\epsilon_k \ \ \ \forall w\in \R^m
\end{equation*}
which, from definition \eqref{eps-subdiff}, is equivalent to $\zk-\ak D_k^{-1}A^T\v \in \partial_{\epsilon_k} f_1^*(A^T\v)$.
Finally, by applying Proposition \ref{pro:zalinescu}, we obtain $A^T\v \in \partial_{\epsilon_k} f_1(\zk-\ak D_k^{-1}A^T\v)$.
Recalling that $\tyk = \zk-\ak D_k^{-1}A^T\v$, which implies $A^Tv = D_k(\zk-\tyk)/\ak$, \eqref{eps-approxk} follows. \endproof

The previous result suggests that for computing $\tyk$ satisfying the assumptions of Corollary \ref{cor:eps-approx} we can use the same iterative approaches described at the end of Section \ref{sec:eta-implementation},
stopping the iterates when
\begin{equation}\label{stop-eps-approx}
\Gskxy\tykl\vl\leq\epsilon_k \ \ \mbox{ and } \thskxy\tykl\xk<0.
\end{equation}

\subsection{Equivalence between $\eta$ and $\epsilon$ approximations}\label{sec:equivalence}

Any $\eta$-approximation $\tyk$ satisfying \eqref{stoprule1} for some $\vl\in\R^m$ is also an $\epsilon$-approximation, where $\epsilon = -\tau\hskxy\tyk\xk$ and $\tau = -1 + 1/\eta$. In fact, in these settings, \eqref{stoprule1} implies $\hskxy\tyk\xk - \Psiskxy\vl\leq -\tau \hskxy\tyk\xk$ and, as shown in Section \ref{sec:eps-implementation}, this means that $\tyk$ is an $\epsilon$-approximation with $\epsilon = -\tau\hskxy\tyk\xk$. Thus, any point computed by an iterative procedure stopped when \eqref{stoprule1} is satisfied, is both an $\eta$- and $\epsilon$- approximation.

\section{Numerical illustration}\label{sec:numerical}

In order to validate the proposed approach, we consider a relevant image restoration problem, whose variational formulation consists in minimizing the sum of a discrepancy functional plus a regularization term. Following the Bayesian paradigm, when the noise affecting the data is of Poisson type, a typical choice for measuring the discrepancy of a given image $\x$ from the observed data $\ve b$ is the following Kullback-Leibler divergence
\begin{equation}\nonumber
KL(\x,\ve b) = \sum_{i=1}^n b_i\log\left(\frac {b_i}{x_i}\right) + x_i - b_i.
\end{equation}
Taking into account also the distortion due to the image acquisition system, which we assume to be modeled through a linear operator $H\in \R^{n\times n}$, and a constant background term $bg$, the data discrepancy is defined as
\begin{equation}\nonumber
f_0(\x) = KL(H\x+bg{\mathds{1}},\ve b),
\end{equation}
where ${\mathds{1}} \in \R^n$ is the vector of all ones.
Moreover, when one wants to preserve edges in the restored image and also the non--negativity of the pixels values, the regularization term can be chosen as
\begin{equation}\nonumber
f_1(x) = \rho \sum_{i=1}^n \|\nabla_i\x\| + \iota_{\R_{\geq 0}^n}(x),
\end{equation}
where $\rho\in \R_{> 0}$ is a regularization parameter multiplying the total variation functional \cite{Rudin-Osher-Fatemi-1992} and $\nabla_i\in \R^{2\times n}$ represents the discrete gradient operator at the pixel $i$. Clearly, the function $f_1(\x)$ has the form \eqref{affine}, with $A=\begin{pmatrix}\nabla_1^T& \cdots&  \nabla_n^T& I\end{pmatrix}^T\in \R^{3n\times n}.$
In this case $\v\in \R^{3n}$ and $g^*$ is the indicator function of the set $B_{0,\rho}^2\times \cdots\times B_{0,\rho}^2\times \R^n_{\leq 0}$, where $B_{0,\rho}^2\subset \R^2$ is the 2-dimensional Euclidean ball centered in 0 with radius $\rho$.

In our experiments we assume that $H$ corresponds to a convolution operator associated to a Gaussian kernel, with reflective boundary conditions, so that the matrix-vector products involving $H$ can be performed via the Discrete Cosine Transform \cite{Hansen-etal-2006}.

We define a set of test problems in the following way: a reference image has been rescaled so that the pixel values lie in a specified range (this is for simulating different noise levels), then it has been blurred by convolution with a Gaussian kernel with standard deviation $\sigma_\mathrm{psf}$ and the background has been added. Finally, Poisson noise has been simulated with the Matlab \verb"imnoise" function, obtaining the noisy blurred image $\ve b$. The details of each test problem are listed in Table \ref{table:1}. The regularization parameter $\rho$ has been manually tuned to obtain a visually satisfactory solution. For each test problem we numerically compute the optimal value $f^*$ by running the considered algorithms for a huge number of iterations, retaining the smallest value found.
\begin{table}
\begin{center}
{\small
\begin{tabular}{lcccccc}
problem         & ref. image & size& range& $\sigma_{\mathrm{psf}}$  & $bg$&$\rho$\\\hline
\verb"cameraman"& Matlab cameraman & $256^2$& $[0,1000]$ & $1.4$& 5&0.0091\\
\verb"micro" & \cite[Figure 8]{Willett-Nowak-2003}& $128^2$&[1,69] & $3.2$&0.5&0.09\\
\verb"phantom"&Shepp-Logan phantom & $256^2$& $[0,1000]$& 1.4& 10   &0.004
\end{tabular}}
\end{center}
\caption{Test problems description}\label{table:1}
\end{table}
\begin{algorithm}[ht]\caption{Variable Metric Inexact Line--search Algorithm (VMILA)}\label{algo:nuSGP}
Choose $0<\alpha_{\min}\leq\alpha_{\max}$, $\mu \geq 1$, $\delta, \beta \in(0,1)$, $\gamma \in [0,1]$, $\eta\in (0,1]$, $\x^{(0)}\in\Omega$. \\
For $k=0,1,2,...$
\begin{itemize}
\item[1.] Choose $\ak\in[\alpha_{\min},\alpha_{\max}]$, $1\leq\mu_k\leq\mu$ and $D_k\in {\mathcal M}_{\mu_k}$;
\item[2.] Compute $\tyk$: compute a dual vector $\vl\in\R^m$ and the corresponding primal vector $\tykl$ such that \eqref{stoprule1} is satisfied, then set $\tyk = \tykl$.
\item[3.] Set $\dk = \tyk-\xk$;
\item[4.] Compute the steplength parameter $\lamk$ with Algorithm \ref{Algo2};
\item[5.] Set $\xkk = \xk+\lamk\dk$.
\end{itemize}
\end{algorithm}
We implement our inexact algorithm, which is summarized in Algorithm \ref{algo:nuSGP}, in Matlab environment with the following settings:

\emph{Step 1, metric selection:} the scaling matrix $D_k$ is chosen mimicking the split-gradient idea \cite{Lanteri-etal-2002}. In particular, at each outer iteration it is defined as the diagonal matrix with positive entries as follows
    \begin{equation}\nonumber
    [D_k]_{ii}= \max\left( \min\left(\frac{\xk_i}{[H^T{\mathds{1}}]_i},\mu_k\right),\frac 1{\mu_k}\right)^{-1}
    \end{equation}
    where $\mu_k = \sqrt{1 + 10^{10}/k^2}$, so that assumption \ref{hypot2prime} is satisfied. We choose a large initial range for the scaling matrix selection to allow more freedom of choice at the first iterates, where the benefits of the scaling matrix are more relevant \cite{Bonettini-etal-2013a}.

\emph{Step 1, steplength selection:} the parameter $\alpha_k$ is chosen by the same strategy used e.g. in \cite{Bonettini-etal-2009,Prato-etal-2012,Prato-etal-2013}, and its value is constrained in the interval $[\alpha_{\min},\alpha_{\max}]$ with $\alpha_{\min} = 10^{-5}$, $\alpha_{\max} = 10^2$.

\emph{Step 2, computation of the approximated proximal point $\tyk$:} we experienced different inner solvers applied on the primal--dual or on the dual formulation of the inner problem. The best performances have been obtained choosing FISTA applied to the dual problem \eqref{dual}, in the variant proposed in \cite{Chambolle-Dossal-2014} which ensures the convergence not only of the objective function values to the optimal one but also of the iterates to the minimum point. In particular, we set $t_l = (l+a-1)/2$, with $a = 2.1$ in \cite[formula (5)]{Chambolle-Dossal-2014}. For brevity, in the following, we report only the results obtained stopping the inner iterates when criterion \eqref{stoprule2} is met, which corresponds to both an $\eta$ and $\epsilon$ approximation (see section \ref{sec:equivalence}). A maximum number of 1500 inner iterations is also imposed. The initial guess of the inner loop at the first outer iterate is the vector of all zeros, while at all successive iterates the inner solver is initialized with the dual solution computed at the previous iterate.

\emph{Other parameters setting:} the line--search parameters $\delta,\beta,\gamma$ have been set respectively equal to $0.5, 10^{-4}, 1$.

All the following results have been obtained on a PC equipped by an Intel Core i7-2620M processor with CPU at 2.70GHz and 8GB of RAM, running Windows 7 OS and MATLAB Version 7 (R2010b).

We investigate first the impact of the inexactness parameter $\eta$ choice on the overall method. In Figure \ref{fig:results1} the relative decrease of the objective function values in the first 500 iterates is reported with respect to both the iteration number (first row) and the computational time, in seconds (second row). It can be observed that a higher precision can accelerate the progress toward the solution, but this usually results in a very large number of inner iterations and, consequently, it is extremely time consuming (for example, for the test problem \verb"cameraman" with $\eta = 10^{-6}, 10^{-2},5\cdot 10^{-1}$ the mean number of inner iterations per outer iteration is 28, 54, 409, respectively). This is typical of inexact algorithms based on the iterative solution of an inner subproblem. We find that a good balance between convergence speed and computational cost is obtained by allowing a relatively large tolerance, corresponding to $\eta = 10^{-6}$. 

As further benchmark, we compare our algorithm to a well established state-of-the-art method, the Chambolle and Pock's method (CP) \cite{Chambolle-Pock-2011}, which, referring to the notations used in their paper, has been implemented setting $G(\x) = \iota_{\R^n_{\geq 0}(x)}$ and $F(K\x) = KL(H\x+bg,\ve b) + \beta \sum_{i=1}^n \|\nabla_i\x\|$, with $K=\begin{pmatrix}H^T, \nabla_1^T,\cdots,\nabla_n^T\end{pmatrix}^T$. In this way the resolvent operator associated to $F^*$ can be computed in closed form. In Figure \ref{fig:results2}, we compare the behaviour of our approach (with $\eta = 10^{-6}$) with CP (2000 iterations) for different choices of its two parameters, $\sigma$ and $\tau$ (once $\tau$ is selected, $\sigma$ is chosen such that $\tau\sigma L^2=1$, where $L = \|K\|$). We can observe that CP is quite sensitive to these parameters, and it is difficult to devise, in general, the more convenient choice, while our approach with the parameters settings described above seems to be always comparable to the best results obtained by CP in terms of objective function decrease with respect to both the iteration number and the computational time.

\def\subFigScale{0.28}
\begin{figure}
\begin{center}
\begin{tabular}{ccc}
\includegraphics[scale=\subFigScale]{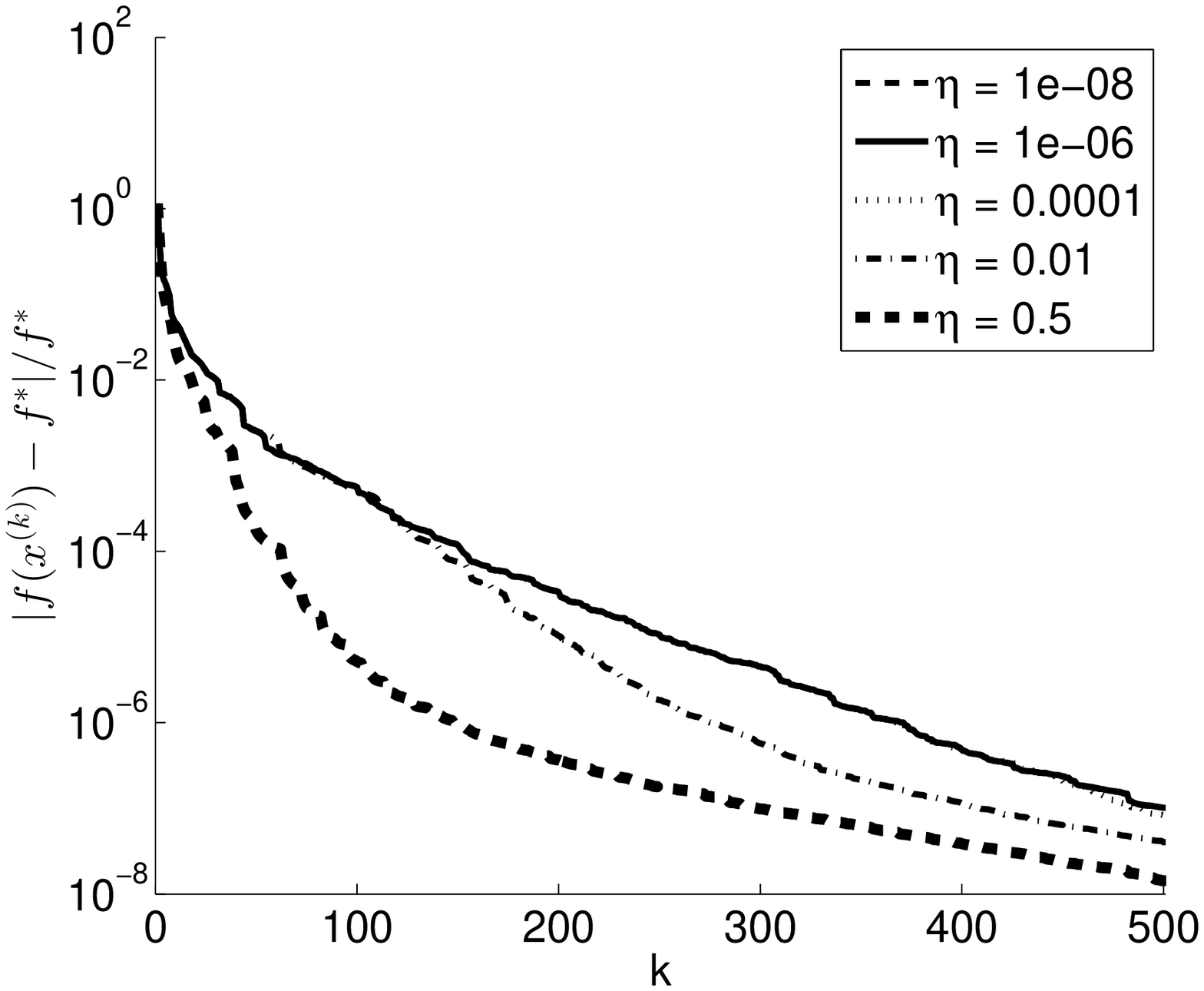}&
\includegraphics[scale=\subFigScale]{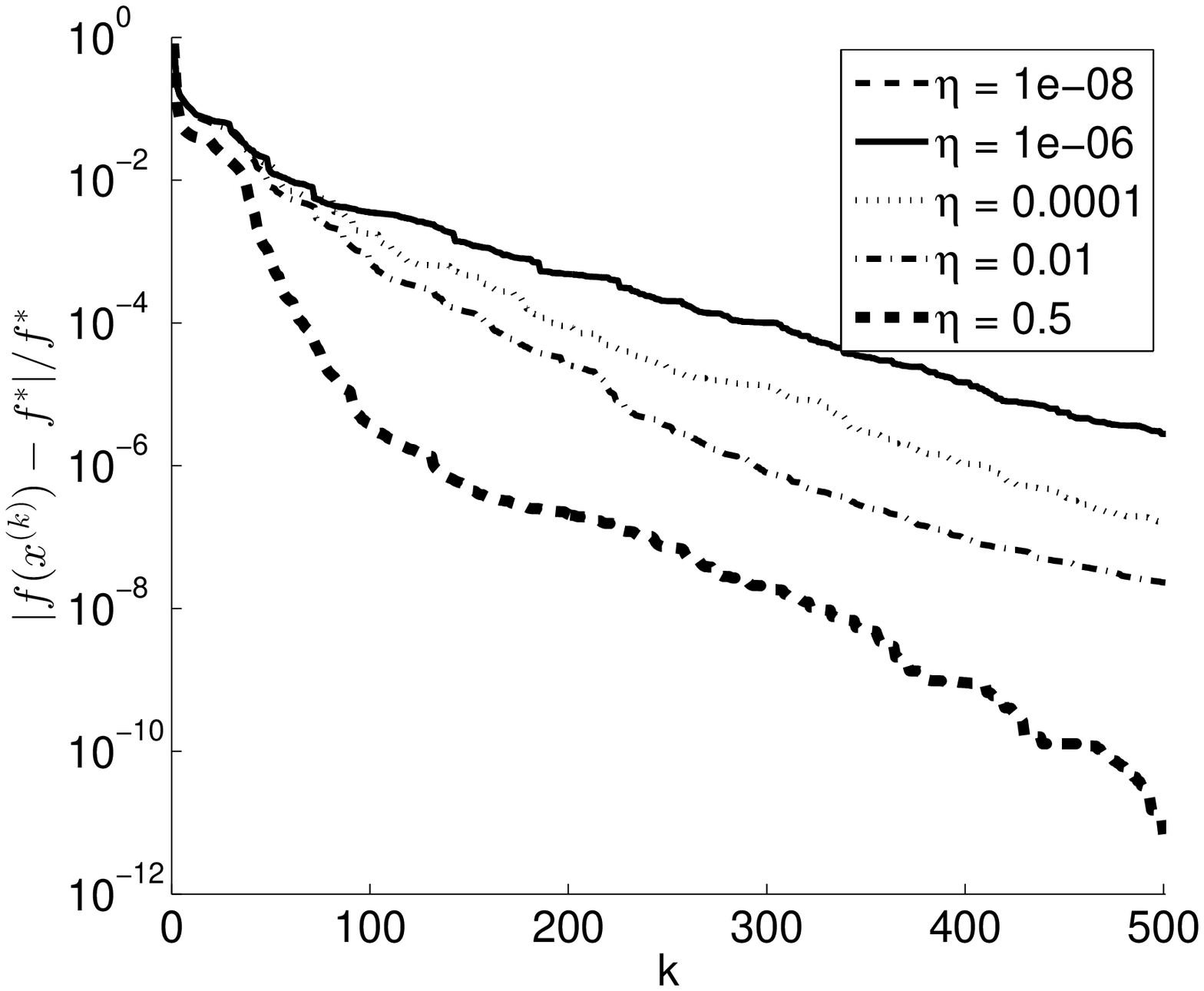}&
\includegraphics[scale=\subFigScale]{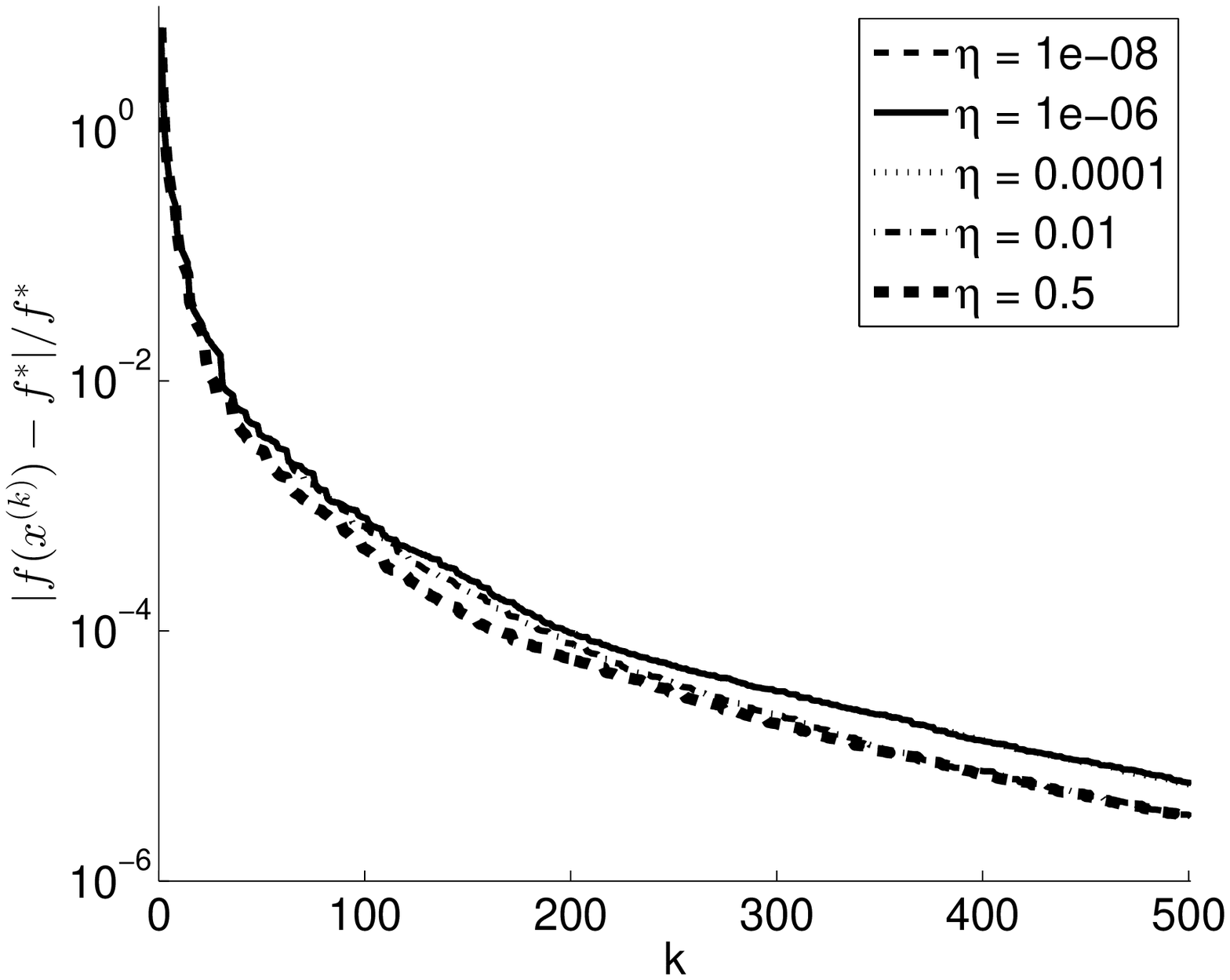}\\
\includegraphics[scale=\subFigScale]{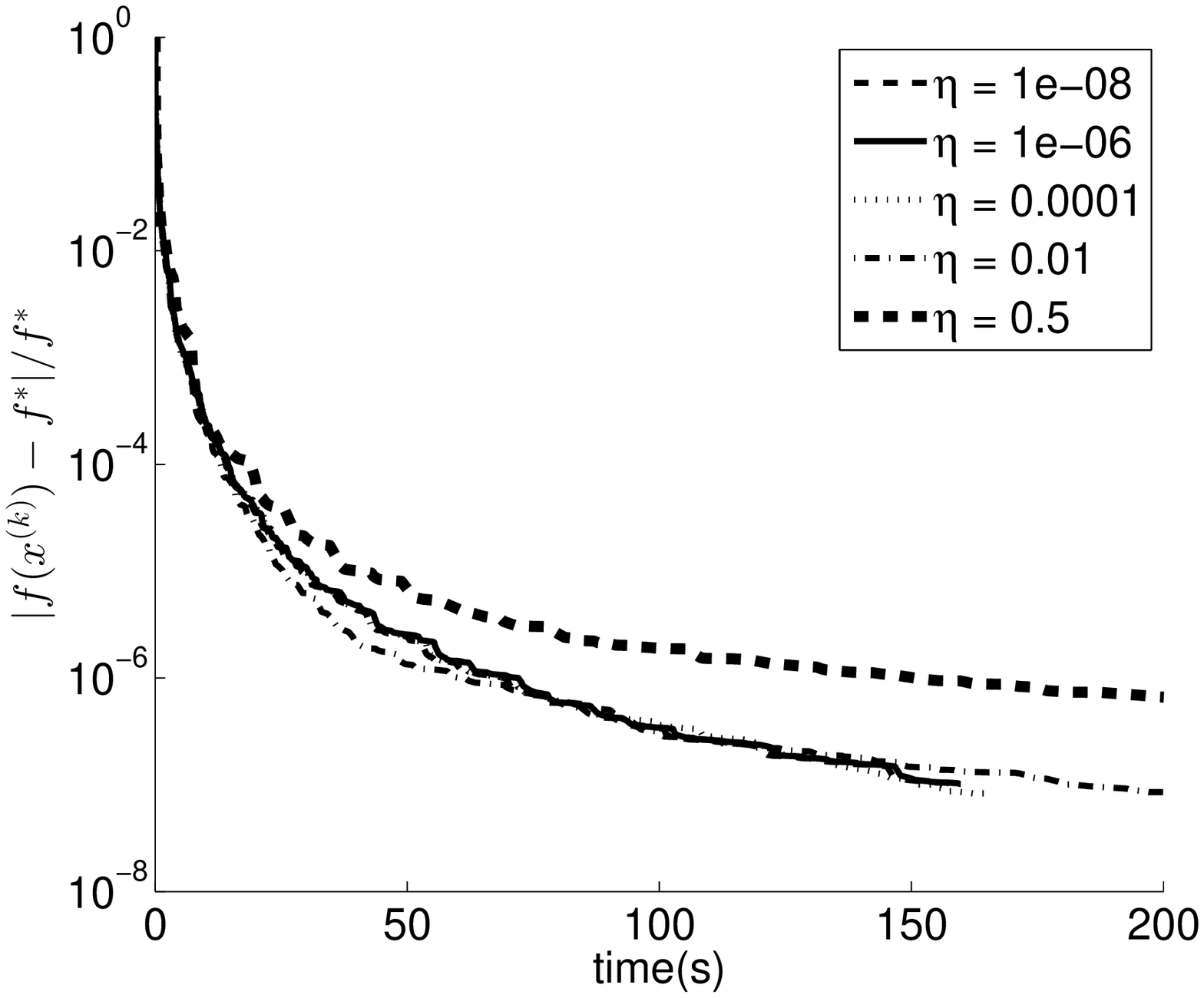}&
\includegraphics[scale=\subFigScale]{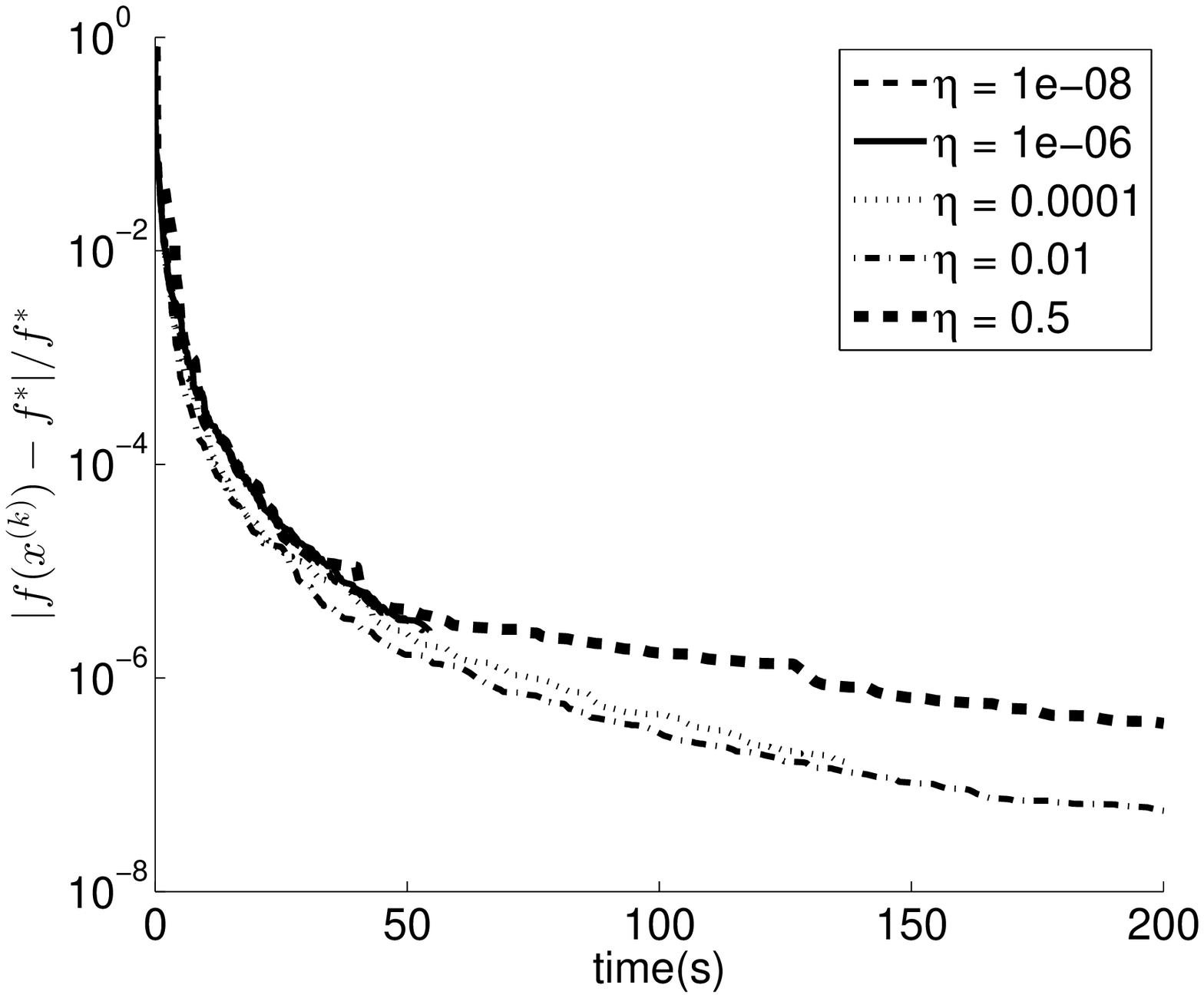}&
\includegraphics[scale=\subFigScale]{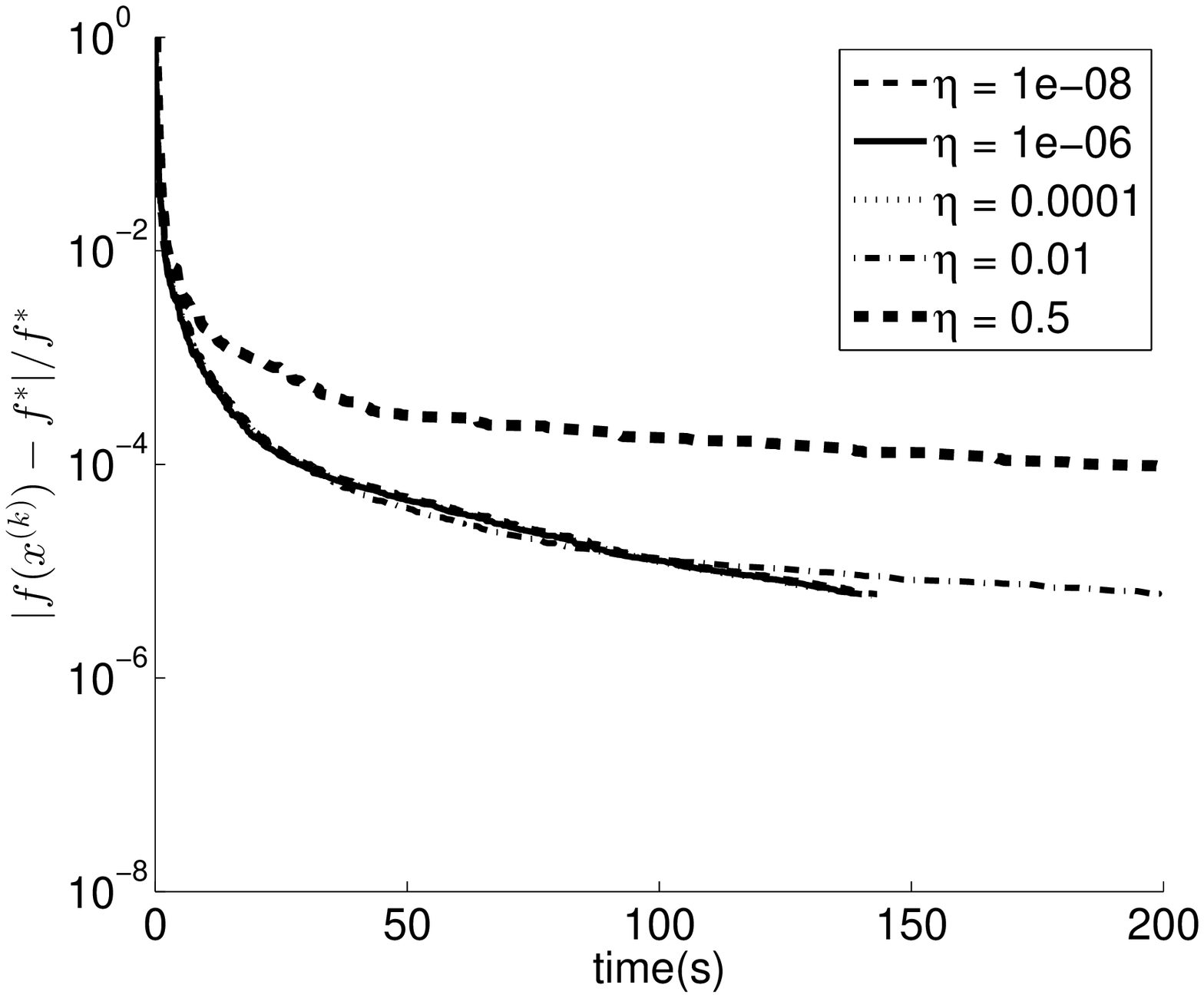}
\end{tabular}\end{center}
\caption{Algorithm \ref{algo:nuSGP} with different choices for $\eta$. Relative decrease of the objective function values with respect to the outer iteration number (top row) and to the computational time (bottom row). Left column: cameraman. Middle column: micro. Right column: phantom.}\label{fig:results1}
\end{figure}

\begin{figure}
\begin{center}
\begin{tabular}{ccc}
\includegraphics[scale=\subFigScale]{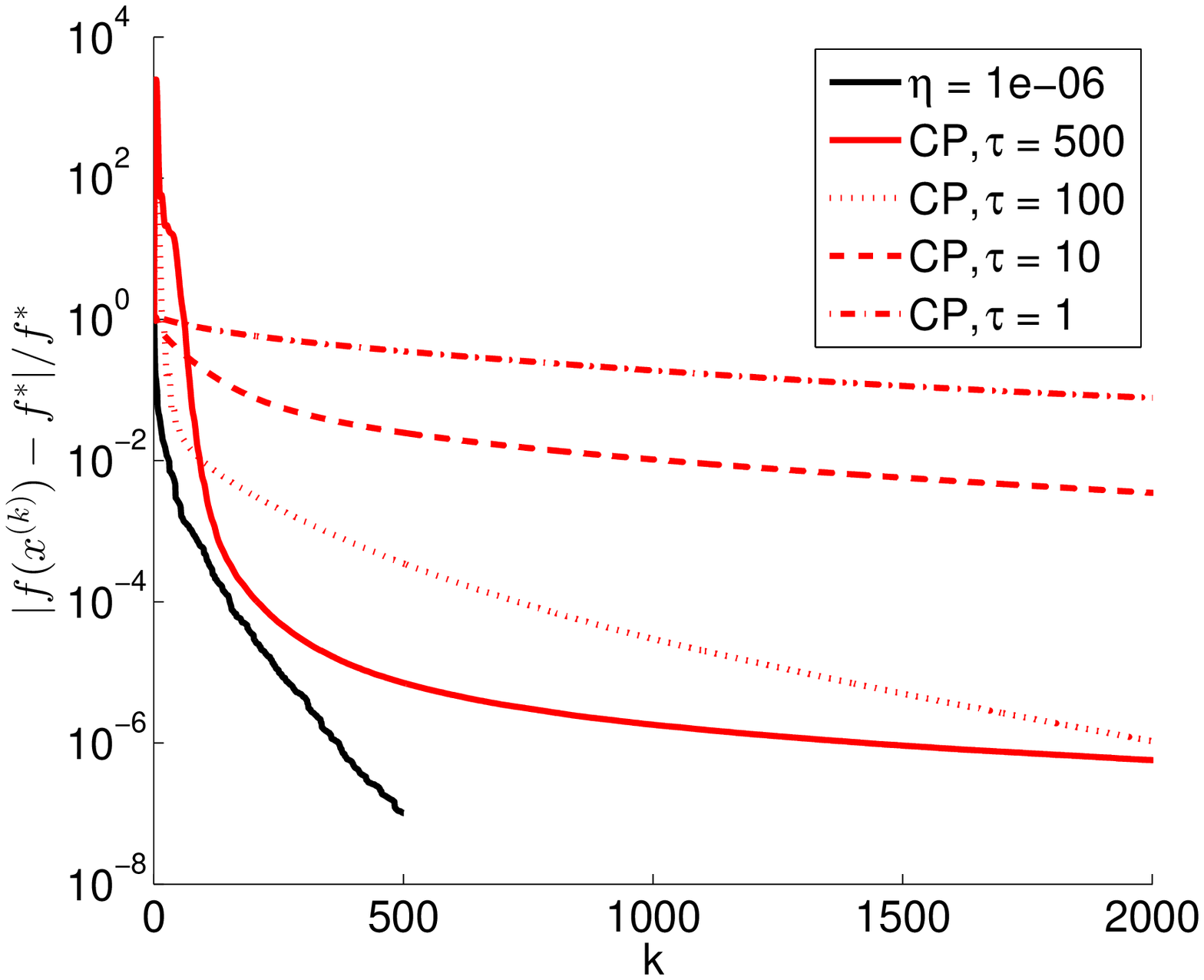}&
\includegraphics[scale=\subFigScale]{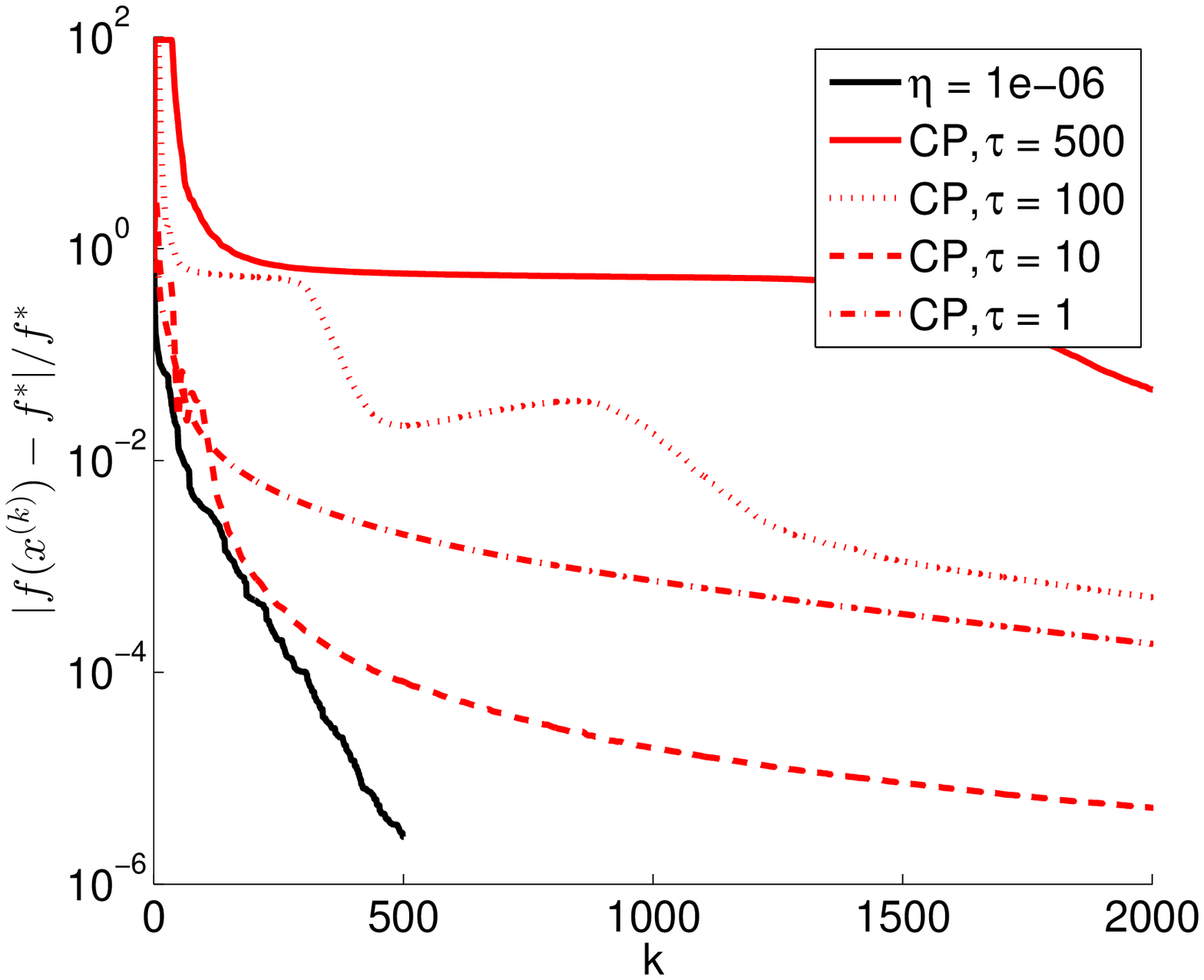}&
\includegraphics[scale=\subFigScale]{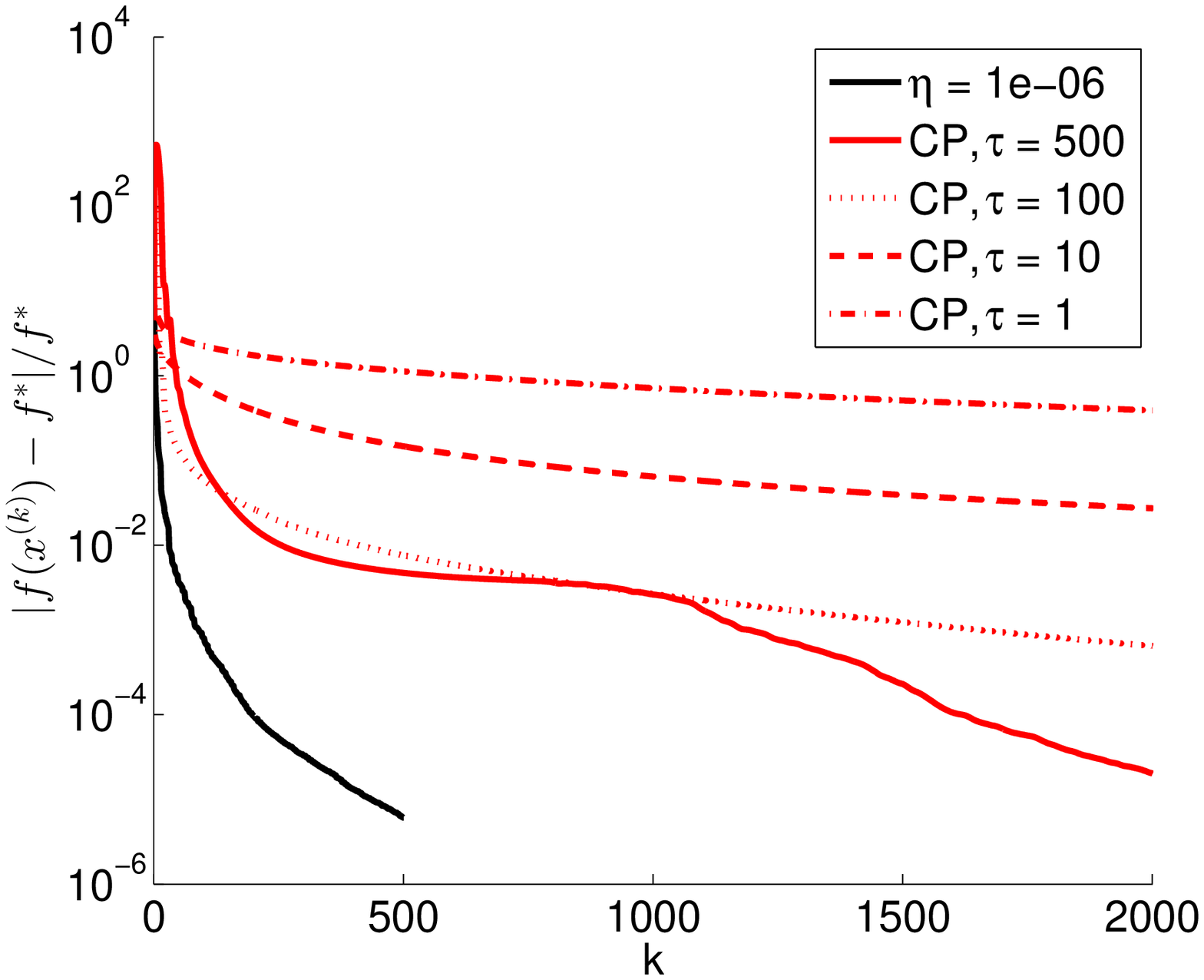}\\
\includegraphics[scale=\subFigScale]{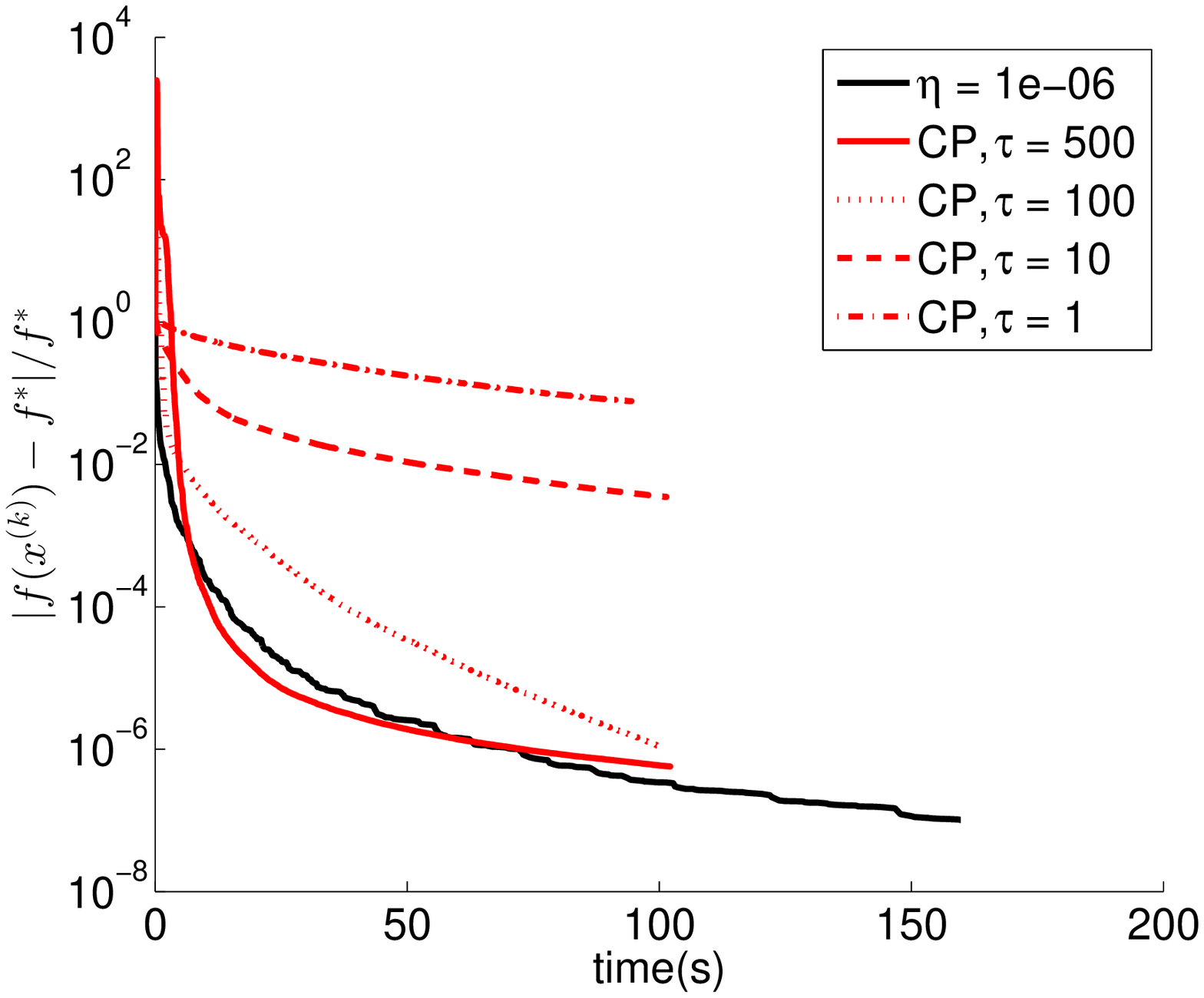}&
\includegraphics[scale=\subFigScale]{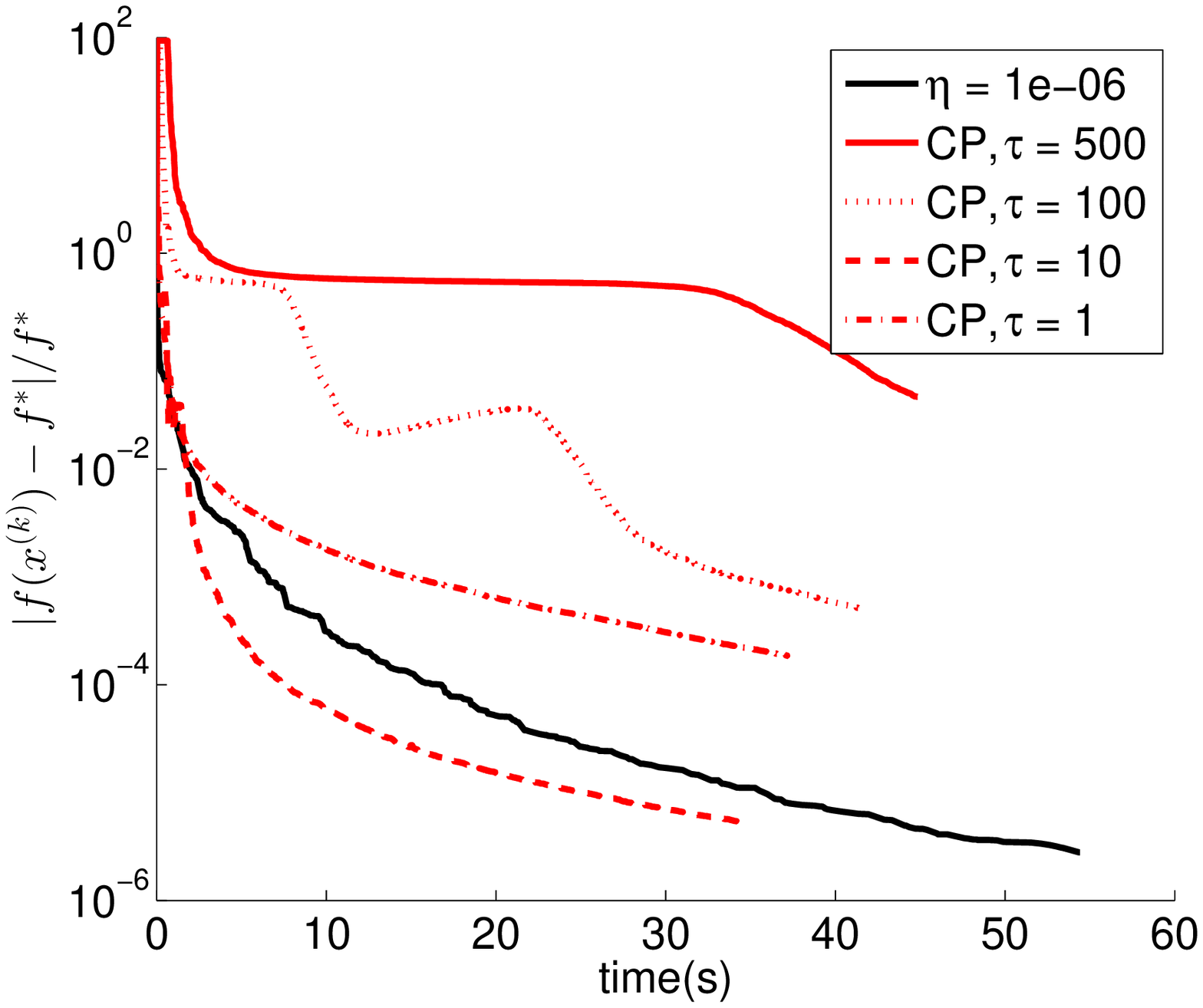}&
\includegraphics[scale=\subFigScale]{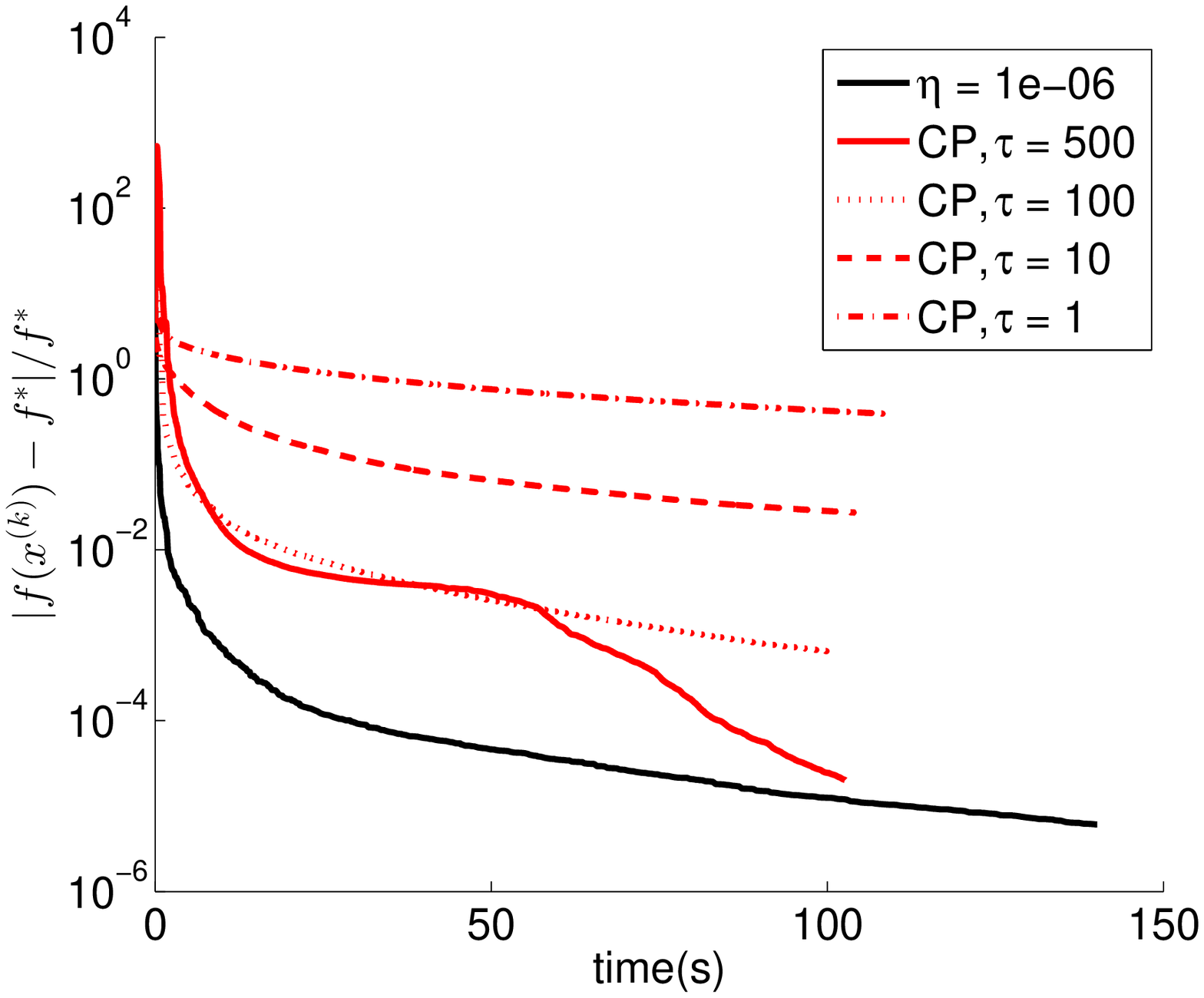}
\end{tabular}\end{center}
\caption{Comparison between Algorithm \ref{algo:nuSGP} ($\eta=10^{-6}$) and the CP algorithm with different choices of its parameters. Relative decrease of the objective function values with respect to the outer iteration number (top row) and to the computational time (bottom row). Left column: cameraman. Middle column: micro. Right column: phantom.}\label{fig:results2}
\end{figure}

\section{Conclusions and future work}\label{sec:conclusions}

In this paper we presented and analyzed an inexact variable metric forward--backward method based on an Armijo--type line--search along a suitable descent direction. The inexactness of the method relies in the possibility of using an approximation of the proximal operator, while the underlying metric may change at each iterations and also non Euclidean metrics are allowed. We performed the convergence analysis of the method, obtaining results in both the nonconvex and convex cases and providing also a convergence rate estimate in the latter one. The main strengths of the method are listed below.
\begin{itemize}
\item The convergence is ensured by a line--search procedure, which does not depend on any user supplied parameter (actually the constants $\gamma,\beta,\delta$ have to be chosen, but the behaviour of the whole algorithm is not sensitive to these choices). On the other side, the ``free'' parameter $\s$ in \eqref{proj} could be exploited to accelerate the convergence speed.
\item The possibility of using at each iterate an approximation of $p(\xk;\hs)$ makes the method well suited for the solution of a wide variety of structured problems.
\item The numerical results on a large scale convex problems shows that the performances of the inexact method are promising and comparable with those of a state-of-the-art method.
\end{itemize}
Future work will be addressed especially to deepen the theoretical and numerical analysis in the nonconvex case, investigating the possibility to obtain convergence results stronger than the ones stated in Theorems \ref{teo-suff-decr} and \ref{thm:eta-approx}, at least for some classes of nonconvex functions (e.g. Kurdyka-{\L}ojasiewicz functions).

\bibliography{biblio_Silvia}
\bibliographystyle{siam}

\end{document}